\documentclass[preprint,12pt,authoryear]{elsarticle}
\usepackage[abs]{overpic}
\newcommand{\cl}{\text{cl}}

\usepackage[ruled]{algorithm2e}
\usepackage[noend]{algpseudocode}

\usepackage{url}

\usepackage{amsmath}
\usepackage{amssymb,mathrsfs}
\usepackage{amsfonts}
\usepackage{nccmath}
\usepackage{graphicx,subfigure}

\usepackage{epstopdf}
\usepackage{hyperref}
\hypersetup{
    colorlinks,
    citecolor=black,
    filecolor=black,
    linkcolor=black,
    urlcolor=black
}

\usepackage[noabbrev,nameinlink]{cleveref}
\usepackage{nameref}

\usepackage{amsthm}
\theoremstyle{plain}
\newtheorem{theorem}{Theorem}[section]
\newtheorem{lemma}[theorem]{Lemma}
\newtheorem{proposition}[theorem]{Proposition}

\theoremstyle{definition}

\newtheorem{remark}{Remark}

\newtheorem{fact}{Fact}[section]

\usepackage{mathtools}
\usepackage[margin=2cm]{geometry}
\usepackage{graphicx}
\graphicspath{ {./images/} }
\usepackage{color}
\usepackage{xcolor}
\usepackage{float}
\usepackage{balance}
\crefformat{equation}{(#1)}
\usepackage{bbm}
\usepackage{bm}
\usepackage{comment}
\usepackage{arcs}

\usepackage{resizegather}

\usepackage{ifthen}
\newcommand{\prob}[2][]{
\ifthenelse { \equal {#2} {} }  %
    { \mathbb{P}_{#1}}   %
    { \mathbb{P}_{#1}\left[#2\right] }   %
}

\usepackage{acronym}
\newacro{lqr}[LQR]{Linear Quadratic Regulators}
\newacro{slqr}[SLQR]{Structured Linear Quadratic Regulators}
\newacro{olqr}[OLQR]{Output-feedback Linear Quadratic Regulators}
\newacro{lqg}[LQG]{Linear Quadratic Gaussian}
\newacro{dare}[DARE]{Discrete-time Algebraic Riccati Equation}
\newacro{ouralgo}[RNPO]{Riemannian Newton-type Policy Optimization}
\newacro{PO}[PO]{Policy Optimization}
\newacro{pg}[PG]{Projected Gradient}

\newcommand{\transpose}{{{\mathsf T}}}
\newcommand{\mK}{{\mathsf{K}}}

\usepackage{textcomp}
\def\BibTeX{{\rm B\kern-.05em{\sc i\kern-.025em b}\kern-.08em
		T\kern-.1667em\lower.7ex\hbox{E}\kern-.125emX}}

\usepackage{enumerate}

\def\bR{{\mathbb{R}}}

\DeclareMathOperator{\grad}{grad}
\DeclareMathOperator{\hess}{Hess}
\DeclareMathOperator{\euchess}{\text{$\overline\hess$}}

\def\lyap{{\,\mathbb{L}}}

\def\fX{\mathfrak{X}}

\DeclareMathOperator{\diff}{\mathrm{d}}

\def\lyaptrace{{Lyapunov-trace}}

\def\Amatrices{\mathbb{R}^{n\times n}}
\def\Bmatrices{\mathbb{R}^{n\times m}}

\def\Lmatrices{\mathbb{R}^{m\times d}}
\def\Kmatrices{\mathbb{R}^{m\times n}}

\def\constraint{{\mathcal{K}}}
\def\stableK{\mathcal{S}}
\def\substableK{\widetilde{\stableK}}

\newcommand{\tensor}[3]{\ensuremath{\left\langle #1, #2 \right \rangle}_{#3}}

\newcommand{\tr}[1]{\ensuremath{\mathrm{tr}\left[ #1 \right]}}

\makeatletter
\newcommand{\algorithmfootnote}[2][\footnotesize]{%
  \let\old@algocf@finish\@algocf@finish%
  \def\@algocf@finish{\old@algocf@finish%
    \leavevmode\rlap{\begin{minipage}{\linewidth}
    #1#2
    \end{minipage}}%
  }%
}
\makeatother

\journal{Encyclopedia of Systems and Control Engineering}

\begin{document}

\begin{frontmatter}

\title{Policy Optimization in Control: Geometry and Algorithmic~Implications}
\author[harvard]{Shahriar Talebi}
\author[uc]{Yang Zheng}
\author[uw]{Spencer Kraisler}
\author[harvard]{Na Li}
\author[uw]{Mehran Mesbahi}

\affiliation[harvard]{organization={Harvard University, School of Engineering and Applied Sciences},%
            addressline={\\150 Western Ave}, 
            city={Boston},
            postcode={02134}, 
            state={MA},
            country={US}}
            
\affiliation[uc]{organization={University of California San Diego, Department of Electrical and Computer Engineering},%
            addressline={\\9500 Gilman Drive}, 
            city={La Jolla},
            postcode={92093}, 
            state={CA},
            country={US}}
            
\affiliation[uw]{organization={University of Washington, Department of Aeronautics and Astronautics},%
            addressline={\\3940 Benton Ln NE}, 
            city={Seattle},
            postcode={98195}, 
            state={WA},
            country={US}}

\begin{abstract}
This survey explores the geometric perspective on policy optimization within the realm of feedback control systems, emphasizing the intrinsic relationship between control design and optimization. By adopting a geometric viewpoint, we aim to provide a nuanced understanding of how various ``complete parameterization''---referring to the policy parameters together with its Riemannian geometry---of control design problems, influence stability and performance of local search %
algorithms. 
The paper is structured to address key themes such as policy parameterization, the topology and geometry of stabilizing policies, and their implications for various (non-convex) dynamic performance measures.
We %
focus on a few iconic control design problems, including the Linear Quadratic Regulator (LQR), Linear Quadratic Gaussian (LQG) control, and $\mathcal{H}_\infty$ control. 
In particular, we first discuss the topology and Riemannian geometry of stabilizing policies, distinguishing between their static and dynamic realizations.
Expanding on this geometric perspective, we then explore structural properties of the aforementioned performance measures and their interplay with the geometry of stabilizing policies in presence of policy constraints; along the way, we address issues such as spurious stationary points, symmetries of dynamic feedback policies, and (non-)smoothness of the corresponding performance measures. We conclude the survey with algorithmic implications of policy optimization in feedback design.

\end{abstract}

\begin{graphicalabstract}
\begin{figure}[h]
\begin{center}
\includegraphics[width= 0.7\textwidth]{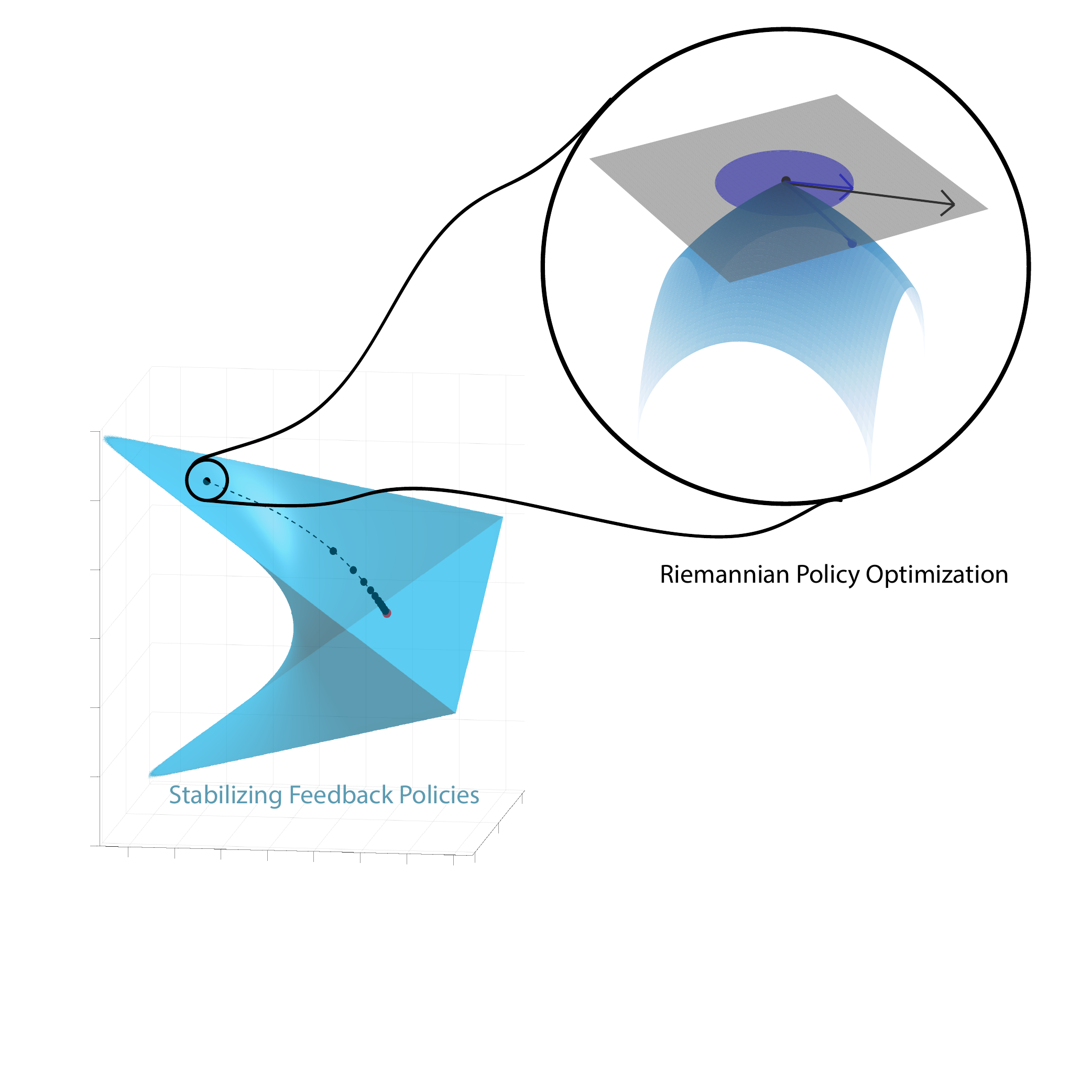}
\end{center}
The article surveys geometric aspects of stabilizing feedback control in the context of policy optimization for iconic synthesis problems, including LQR, LQG, and ${\cal H}_{\infty}$.
\end{figure}
\end{graphicalabstract}

\begin{highlights}
\item Policy optimization (PO) provides a unified perspective on feedback design subject to stabilization constraint,
\item PO also provides a bridge between control theory and data-driven design techniques such as model-free reinforcement learning,
\item Geometric properties of the set of stabilizing feedback gains and how their interact with various performance measures become of out most importance while designing PO-based algorithms,
\item PO highlights design and algorithmic challenges in control engineering beyond classic paradigms involving quadratic costs and linear dynamics.
\end{highlights}

\end{frontmatter}
\tableofcontents
\newpage
\section{Introduction}
\label{intro}
Optimization and control have had a rich symbiotic relationship since their inception. This is not surprising as the current dominate perspective on control design highlights: the design process for a dynamic system involves formalizing the notion of ``best'' with respect to selectable system parameters\footnote{Whatever the ``best'' qualifier implies.}--followed by devising algorithms that optimize the design objective with respect to these parameters. 
The interactions between these two disciplines has highlighted a crucial aspect of control design: identifying design objectives that, with respect to a given parameterization of control system (or more generally, parameters in the overall design architecture), facilitate algorithmic developments while maintaining sound engineering judgment.
For example, it is common to see a family of solution strategies for a given control design problem, including open vs. closed-loop, approaches based on variational~\citep{Liberzon2011-en} or dynamic programming~\citep{Bertsekas2012-oq}, and constructs such as co-state, value, and, of course, {\em policy}~\citep{Bertsekas2011-hi,Bertsekas2017-fi}. Each of these formalisms offer insights into control design, that although in principle can be formalized as an optimization problem, but reflect distinct intricacies in designing systems that evolve over time, interact with their environments, and have memory. For example, system theoretic properties such as stability, minimality, and robustness, significantly ``spice up" not only the control design problem formulation but also the adopted strategies for their solution~\citep{Zabczyk_undated-jw,Sontag2013-hk}.

A powerful abstraction in this plethora of design techniques is that of \textit{policy}, mapping what the system has observably done 
to what can influence its subsequent behavior. Not only does the existence of an ``optimal'' policy reduce the complexity of control implementation (as the control input is a trajectory) at the expense of realizing feedback, but as it turns out, it also addresses one of the key features of control design, namely robustness~\citep{Skogestad2005-wi,Doyle2013-vh}. This survey aims to capture the {\em geometry} of policy optimization for a few iconic design techniques in control (including stabilization, linear quadratic control, and $\mathcal{H}_\infty$ robust control), capturing distinct facets of adopting such a perspective for feedback design. However, as we undertake such an endeavor, it is important to comment on the timeliness of this approach, as well as how it connects--and is distinct--from previous works, particularly in relation to its reference to ``geometry''.

Historically, the ``geometric'' qualifier has been used in systems and control to shed light on more subtle aspects of system design. In linear geometric control, subspace geometry is used to delineate the coordinate-independence of system theoretic constructs such as controllability and observability~\citep{wonham_w_murray_linear_1985,Trentelman2001-ck}. Geometric nonlinear control, on the other hand,  characterizes notions such as controllability for nonlinear systems, by viewing their evolution in terms of vector fields, and then brings forth a differential geometric formalism for their analysis, e.g., differentiable manifold, distributions, and Lie algebras~\citep{Isidori2013-dj,Nijmeijer1990-xr}; the aim is to free control theory (analysis and synthesis) from the confinements of linearity and linear algebra via a coordinate-free analysis of dynamical systems. Using this geometric vista, linear maps are uplifted to diffeomorphisms, while notions
such as controllability matrix are revealed as Lie brackets~\citep{Brockett2014-zg}. As these two examples demonstrate, ``geometry'' is often used to hint at coordinate-independence of concepts, similar to how finite dimensional vector spaces are related to linear algebra. Other notable works in adopting a geometric perspective in systems and control theory, particularly in relation to realization theory and system identification include~\citep{beckmann_invariants_1976,byrnes_fine_1980,tannenbaum_invariance_1981}.

By including “geometry” in the title of this survey, we deliberately mean
to promote adopting a similar geometric perspective as the aforementioned works, but for the space of feedback control policies rather than system models or their trajectories, which highlights the importance of how one characterizes notions such as distance and direction for these policies.
Such a perspective not only complements other features of this space (in addition to its topological, analytic, or algebraic structures) but more importantly, has direct consequences for devising algorithms for the corresponding optimization problems.
This geometric perspective and its algorithmic implications have also been adopted in neighboring decision sciences, e.g., statistical learning and Markov Decision Processes (MDPs). Notable in the landscape of such geometric techniques, we mention the notion of natural gradients, where the geometry of the underlying model, singled out in the design objective, is systematically used to synthesize algorithms that behave invariant under certain re-parameterizations. By invariance, we mean embedding the underlying model with a notion of distance that is preserved under certain mappings, e.g., Fisher metric in statistical learning~\citep{amari_natural_1998,Amari2016-do}. Closely related to the present survey is adopting the theory of natural gradients for MDPs as first proposed in~\citep{Kakade2002-fc}, where the Hessian geometry induced by the entropy plays a central role in the
design and convergence analysis of the corresponding algorithms; also see~\citep{Agarwal2021-ss}. 

These geometric insights have a number of algorithmic and system-theoretic consequences. For example, as we will see, improving the policy by taking steps in the direction of its (negative) gradient proves to be an effective means of obtaining optimal policies, i.e., first-order policy updates. The ``zeroth" order version of the above scheme, on the other hand, leads to using function evaluations to approximate the corresponding gradients from data--say, when such evaluations can be obtained from an oracle that can return approximate values of the cost, closely relates to reinforcement learning setup.
Key questions in this data-driven realization of geometric first-order methods are with how many function evaluations (and with what accuracy) and over how many iterations, an accurate estimate of the optimal policy can be obtained~\citep{Fazel2018-pv,Malik2019-bs,Mohammadi2021-vs,Mohammadi2020-iv}. Furthermore, these concepts resonate with optimal estimation problems due to the profound duality relation between control and estimation \citep{talebi_data-driven_2023}.

The survey is structured as follows. In \S\ref{parameterization} we make the parameterization theme on control design alluded to above more explicit. While this perspective provides a direct approach to formalize policy optimization for control, it also underscores how the constraints imposed by system theoretic notions such as stabilizability make the feasible set of the optimization non-trivial. This is first more pursued for static feedback, followed by that of dynamic feedback policies in \S\ref{section:stablizing-policies}. We then turn our attention to how distinct design performance measures interact with the feasible set of feedback policies in \S\ref{section:LQ-performance}.
Algorithmic implications and data-driven realizations of the geometric perspective on policy optimization are then examined under \S\ref{data-driven}. In \S\ref{summary}, we provide a summary of the key points put forth by this survey as well as our outlook on the future work; \S\ref{notes-commentary} provides commentary on references with contributions reflected in this survey.

\section{Policy Optimization in Control: the Role of Parameterization} \label{parameterization}

Let us consider a discrete-time dynamical system:
\begin{equation}\label{eq:system-dynamics}
         x_{t+1} = f(x_t,u_t,w_t), \quad y_t = h(x_t,w_t),
\end{equation} 
where $x_t \in \mathbb{R}^n$ is the system state, $u_t \in \mathbb{R}^m$ is the control input, and $y_t \in \mathbb{R}^p$ is the system measurement. We refer to $w_t$ as the exogenous signal,
representing unmodeled dynamics, stochastic noise, or disturbances.
At each time step $t$, we use $l(x_t,u_t)$ to denote the stage cost as a function of the current state and input. The goal of infinite-horizon optimal control is to choose $\mathbf{u}=(u_0,u_1,\ldots)$ to minimize an accumulated cost over the infinite time horizon. 
More formally, we define the $T$-stage accumulated cost as

\begin{align}\label{eq:ave-inf-hor-cost}
    J_{T}(\mathbf{u}, \mathbf{w}, x_0) \coloneqq  \sum_{t=0}^T l(x_t,u_t),    
\end{align} where $\mathbf{w}=(w_0,w_1,\ldots)$. The cost (\ref{eq:ave-inf-hor-cost}) highlights that the implicit state trajectory $\mathbf{x}$ depends on the control input $\mathbf{u}$, the exogenous input $\mathbf{w}$, and initial state $x_0$.
Often $\mathbf{w}$ and $x_0$ are modeled either stochastically or deterministically, and the performance of the closed-loop system is measured based on the corresponding average or worst case performance. 
An example of such infinite-horizon control performance is
\begin{equation}\label{eq:single-param-ave-inf-hor-cost}
    J(\mathbf{u}) := \lim_{T \to \infty} \mathbb{E}_{\mathbf{w},x_0} \;\frac{1}{T} J_T(\mathbf{u},\mathbf{w}, x_0),
\end{equation}
which presumes stochastic exogenous input $\mathbf{w}$ and initial state $x_0$ 
(by expectation with respect to statistical properties of $\mathbf{w}$ and $x_0$); as in the formulation of LQR and LQG costs. 
On the other hand, when the exogenous inputs are adversarial, one may replace the expectation with respect to $\mathbf{w}$ with the worst-case performance assuming bounded energy for $\mathbf{w}$; as in the formulation of $\mathcal{H}_\infty$ cost.

\subsection{Policy Parameterization for Closed-loop Optimal Control}\label{subsect:policy-param-costs}

Instead of optimizing over the input sequence $\mathbf{u}$ in (\ref{eq:single-param-ave-inf-hor-cost}), which we refer to as the \textit{optimal open-loop control design}, we consider instead
optimizing over a class of \textit{feedback policies} that act
on the system history $\mathcal{H}_t := (u_{0:t-1}, y_{0:t})$, where $u_{0:t-1} := (u_0,u_1,\ldots,u_{t-1})$ and similar for $y_{0:t}$. As such a \textit{feedback} or \textit{closed-loop policy} at time $t$, denoted by $\pi_t: \mathcal{H}_t \mapsto u_t$, is a measurable function that maps the system history $\mathcal{H}_t$ at time $t$ to a control input $u_t$. We can alternatively define $\pi_t(\mathcal{H}_t)$ to be a distribution and set $u_t \sim \pi_t(\mathcal{H}_t)$.\footnote{In this case, we would have to slightly augment (\ref{eq:single-param-ave-inf-hor-cost}).} We will call $\pi=(\pi_0,\pi_1,\ldots)$ the feedback policy;
for brevity, we will often write $\pi(\mathcal{H}_t):=\pi_t(\mathcal{H}_t)$. 
Let $\Pi$ be the infinite-dimensional vector space of all such feedback policies $\pi(\cdot)$. 
Then, the optimal closed-loop policy problem reads as 
\begin{subequations}\label{eq:closed-loop-policy-prob}
    \begin{align}
        \min_{\pi \in \Pi} \quad & \ J(\mathbf{u}) \\
        \text{subject to} \quad &(\ref{eq:system-dynamics}),\;  
        u_t = \pi(\mathcal{H}_t). %
    \end{align}
\end{subequations} 

Since it is non-trivial to optimize directly over this space $\Pi$, the so-called ``policy parameterization'' approach is to parameterize (a subset of) $\Pi$ with some $d$-dimensional set of parameters $\theta \in \Theta \subset \mathbb{R}^{d}$. That way, our parameterized family of policies $\{\pi_\theta(\cdot)\}_{\theta \in \Theta} \subset \Pi$ is finite-dimensional. Note that policy parameterization is rather flexible, as it can represent linear dynamical systems, polynomials, kernels, or neural networks. In some important cases, one can restrict the class of policies under consideration without loss of generality. For instance, it is known that static/dynamic ``linear policies'' are sufficient for optimal and robust control problems posed for linear time-invariant (LTI) systems~\citep{zhou1996robust}. With policy parameterization, \eqref{eq:closed-loop-policy-prob} is reduced to the optimal closed-loop parameterized policy problem:

\begin{equation} \label{eq:closed-loop-param-policy-prob}
    \min_{\theta \in \Theta}  J(\theta) := J(\pi_{\theta}(\mathcal{H})), 
\end{equation} 
where $\pi_{\theta}(\mathcal{H})$ denotes the input signal obtained as $\mathbf{u} = (\pi_\theta(\mathcal{H}_0), \pi_\theta(\mathcal{H}_1), \dots)$. If $\pi^*$ solves (\ref{eq:closed-loop-policy-prob}) and $\theta^*$ solves (\ref{eq:closed-loop-param-policy-prob}) then it should be remarked $J(\pi^*) \leq J(\pi_{\theta^*})$, with equality if and only if the parameterization is ``rich enough,'' e.g., static and dynamic linear feedback policies for LQR and LQG problems, respectively.
Further, we can explicitly incorporate a policy constraint on $\Theta$ that represents closed-loop stability or an information structure required for control synthesis.

Conceptually, it appears simple and flexible to use local search algorithms, such as policy gradient or its variants, to seek the ``best'' policy in \eqref{eq:closed-loop-param-policy-prob}. Once the corresponding cost in \eqref{eq:closed-loop-param-policy-prob} can be estimated from sampled trajectories, this policy optimization setup indeed is very amenable for data-driven design paradigms such as reinforcement learning \citep{recht2019tour}. 
One goal of this article is to highlight some \textit{rich} and \textit{intriguing}, geometry in policy optimization \eqref{eq:closed-loop-param-policy-prob}, including nonconvexity, potentially disconnected feasible domain $\Theta$, spurious stationary points of the cost $J$, symmetry, and smooth (Riemannian) manifold structures. Our focus will be on classic control tasks for LTI systems including LQR, LQG, and  $\mathcal{H}_\infty$ robust control. These geometrical understandings will often provide insights for designing principled local search algorithms to solve \eqref{eq:closed-loop-param-policy-prob}, such as policy gradient methods with global/local convergence guarantees \citep{fazel2018global,bu_lqr_2019,mohammadi2021convergence,hu2023toward}. Inspired by the rich geometry in \eqref{eq:closed-loop-param-policy-prob}, we will further emphasize that a ``complete policy parameterization'' should come with an associated metric, capturing the inherent geometry, which may help improve the problem's conditioning \citep{talebi_riemannian_2022,kraisler2024output}.

\subsection{Policy Optimization for Iconic Optimal and Robust Control Problems}

\subsubsection{LQR under Static State-feedback Policies}

For the LQR problem, the system dynamics \eqref{eq:system-dynamics} read as 
\begin{equation} \label{eq:state-evolution}
    {x}_{t+1}= Ax_t+Bu_t+ w_t, \qquad y_t = x_t,
\end{equation} 
where the process noise is white and Gaussian $w_t \sim \mathcal{N}(0, \Sigma)$ for some $\Sigma \succ 0$.
Since the state is directly observed, the system history takes the form $\mathcal{H}_t = (u_{0:t-1}, x_{0:t})$. Given the performance weighting matrices $Q \succeq 0, R \succ 0$, (\ref{eq:closed-loop-policy-prob}) reduces to the optimal LQR problem:
\begin{subequations}\label{eq:LQR-noise}
    \begin{align}
    \min_{\pi \in \Pi} \quad &J_{\texttt{LQR}}(\mathbf{u}) \coloneqq \lim_{T \to \infty} \mathbb{E}_{\mathbf{w},x_0} \; \frac{1}{T}J_T(\mathbf{u},\mathbf{w}, x_0) \\
    \text{subject to}\quad &(\ref{eq:state-evolution}), \; u_t = \pi(\mathcal{H}_t)
    \end{align} 
\end{subequations} with the quadratic stage cost $l(x_t,u_t)=\frac{1}{2}(x_t^TQx_t + u_t^TRu_t)$. 
Note that the cost in \eqref{eq:LQR-noise} will be oblivious to any bounded initial condition $x_0$ as long as the policy is stabilizing, and the Gaussian assumption on noise can be lifted.

Alternatively, when the dynamics is noiseless but the initial condition is uncertain, we might be interested in minimizing the following objective instead 
\begin{subequations}\label{eq:LQR-noiseless}
    \begin{align}
        \min_{\pi \in \Pi} \quad &J_{\texttt{LQR}}(\mathbf{u}) \coloneqq \lim_{T \to \infty} \mathbb{E}_{x_0}\;  J_T(\mathbf{u}, \mathbf{0}, x_0) \\
        \text{subject to} \quad &(\ref{eq:state-evolution}), \; u_t = \pi(\mathcal{H}_t)
    \end{align}
\end{subequations} where $x_0$ is drawn from a distribution with covariance $\Sigma \succ 0$. We will see later that both of the problems in (\ref{eq:LQR-noise}) and (\ref{eq:LQR-noiseless}) essentially amount to the same policy optimization problem. 

\begin{figure}[t]
    \centering
    \includegraphics[width=.4\textwidth]{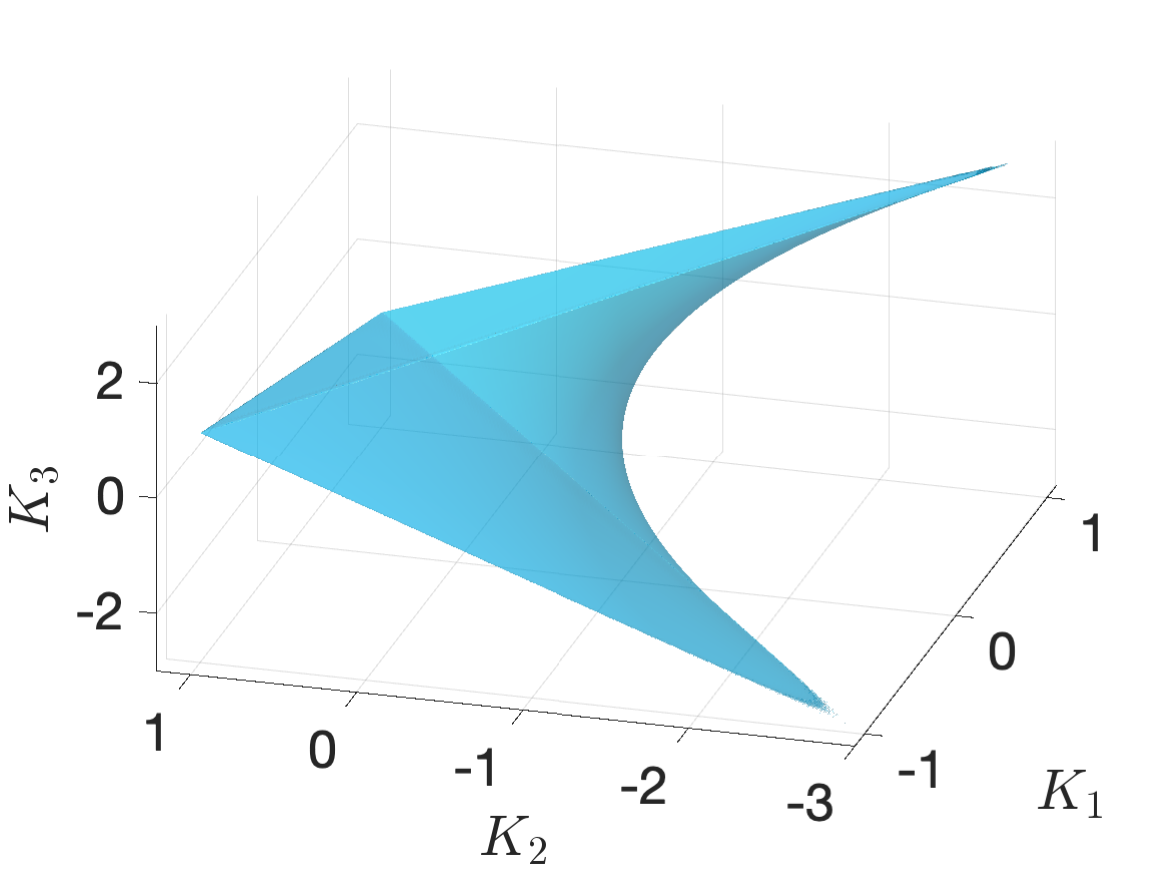}
    \caption{The non-convex set of stabilizing static state-feedback polices $\stableK$ for $A = \begin{bmatrix} 0 & 1 & 0 \\ 0 & 0 & 1 \\ 0 & 0 & 0\end{bmatrix}$ and 
$B = \begin{bmatrix} 0 \\ 0 \\ 1\end{bmatrix}$.
}
    \label{fig:LQR-simple}
\end{figure}

A remarkable property of the closed-loop optimal LQR problem is that the optimal policy $\pi^*$, when exists, is linear in the states and depends only on the current states. That is, $\pi^*(\mathcal{H}_t) = \pi^*(x_t) = K x_t$ for some $K \in \mathbb{R}^{m \times n}$ \citep{anderson_optimal_2007}.
Inspired by this property, for both  \eqref{eq:LQR-noise} and \eqref{eq:LQR-noiseless}, we can parameterize the LQR policy as a linear mapping $x_t \mapsto u_t = Kx_t$ and referred to as policy $K$ for simplicity. Then, the set of static stabilizing state-feedback policies is  
\begin{equation} \label{eq:static-stabilizing-policies}
    \Theta := \mathcal{S} = \{K \in \mathbb{R}^{m \times n}: \rho(A+BK) < 1\}.
\end{equation} where $\rho(\cdot)$ denotes the spectral radius of a square matrix. Following (\ref{eq:closed-loop-param-policy-prob}), we can define $J_{\texttt{LQR}}$ over $\mathcal{S}$, where $J_{\texttt{LQR}}(K)$ corresponds to the objective in either~\eqref{eq:LQR-noise} or \eqref{eq:LQR-noiseless}. Our parameterized policy family will be $\{\pi_K\}_{K \in \mathcal{S}}$ with $\pi_K(\mathcal{H}_t) := K x_t$.
Lastly, the optimal LQR policy problem becomes 
\begin{equation}
\begin{aligned} \label{eq:LQR-policy-optimization}
    \min_{K} \quad & J_{\texttt{LQR}}(K) \\
    \mathrm{subject~to}\quad &  K \in \mathcal{S}. 
\end{aligned} 
\end{equation}
For well-posedness of the problem, we assume that the pair $(A,B)$ is stabilizable so $\stableK$ is non-empty; see \Cref{fig:LQR-simple} for a numerical example.

\subsubsection{LQG under Dynamic Output-feedback Policies}\label{subsubsect:LQG}

We here specify the general closed-loop optimal policy problem \eqref{eq:closed-loop-policy-prob} in the context of LTI systems with partial observation:
\begin{equation} \label{eq:dynamics_discrete}
    \begin{aligned}
        x_{t+1} = A x_t + B u_t + w_t, \qquad  
        y_t = Cx_t + v_t, 
    \end{aligned}
\end{equation}
where $x_t \in \mathbb{R}^n, u_t \in \mathbb{R}^m, y_t \in \mathbb{R}^p$ are the system state, input, and output measurement at time $t$, and $w_t \sim \mathcal{N}(0, W), v_t \sim \mathcal{N}(0, V)$ are Gaussian process and measurement noise signals, respectively.  It is assumed that the covariance matrices satisfy $W \succeq 0, V \succ 0$.
Given performance weight matrices $Q \succeq 0, R \succ 0$, the optimal LQG problem becomes,
\begin{subequations}\label{LQG-policy-problem}
    \begin{align}
        \min_{\pi \in \Pi} \quad &J_{\texttt{LQG}}(\mathbf{u}) \coloneqq \lim_{T \to \infty} \mathbb{E}_{\mathbf{w},x_0} \; \frac{1}{T} J_T(\mathbf{u},\mathbf{w}, x_0) \\ \label{eq:LQG-cost}
        \text{subject~to} \quad &(\ref{eq:dynamics_discrete}), \; u_t = \pi(\mathcal{H}_t)
    \end{align}
\end{subequations} where $\mathbf{w} := ((w_t,v_t))_{t=0}^\infty$. It should be noted that the difference between the optimal closed-loop policy problems (\ref{LQG-policy-problem}) and (\ref{eq:LQR-noise}) is the policies have to to be output-feedback with $\mathcal{H}_t = (u_{0:t-1},y_{0:t-1})$ versus state-feedback with $\mathcal{H}_t =(u_{0:t-1},x_{0:t})$, respectively. 

Next, we construct a family of policies, referred to as \textit{dynamic policies}, parameterized by an LTI system,
\begin{equation} \label{eq:controller_discrete}
    \begin{aligned}
        \xi_{t+1} = A_{\mK} \xi_t + B_{\mK} y_t, \quad 
        u_t = C_{\mK} \xi_{t}, 
    \end{aligned}
\end{equation}
where $\xi_t \in \mathbb{R}^q$ is the controller's internal state at time $t$, and policy parameters $(A_{\mK},B_{\mK},C_{\mK}) \in \mathbb{R}^{q \times q} \times \mathbb{R}^{q \times p} \times \mathbb{R}^{m \times q}$ specify the policy parameters. If $q = n$, we refer to \eqref{eq:controller_discrete} as a \textit{full-order} dynamic policy; if $q < n$, it is called a \textit{reduced-order} dynamic policy. In this survey, we will focus on the case $q = n$ since it is known that full-order dynamic policy parameterization is rich enough; i.e. contains a globally optimal solution to the closed-loop optimal policy problem \eqref{LQG-policy-problem} \citep{zhou1996robust}. 
Therefore, combining~\eqref{eq:controller_discrete} with~\eqref{eq:dynamics_discrete} leads to the augmented closed-loop system,
\begin{equation}\label{eq:closed_loop_system}
\begin{aligned}
   \begin{bmatrix} x_{t+1} \\ \xi_{t+1} \end{bmatrix} = \begin{bmatrix}
    A & BC_{\mK} \\
    B_{\mK}C & A_{\mK}
    \end{bmatrix} \begin{bmatrix} x_t \\ \xi_t \end{bmatrix}  +  \begin{bmatrix} I & 0 \\ 0 & B_{\mK}  \end{bmatrix}\begin{bmatrix} w_t \\ v_t \end{bmatrix}, \quad 
    \begin{bmatrix} y_t \\ u_t \end{bmatrix} = \begin{bmatrix} C & 0 \\ 0& C_{\mK} \end{bmatrix} \begin{bmatrix} x_t \\ \xi_t \end{bmatrix} + \begin{bmatrix} v_t \\0 \end{bmatrix}.
\end{aligned}
\end{equation}
The set of stabilizing controllers with order $q \in \mathbb{N}$ is now defined as, %
\begin{equation} \label{eq:internallystabilizing}
    \Theta := \mathcal{C}_{q} = \left\{
    \left.\mK=\begin{bmatrix}
    0_{m\times p} & C_{\mK} \\
    B_{\mK} & A_{\mK}
    \end{bmatrix}
    \in \mathbb{R}^{(m+q) \times (p+q)} \right|\;  \rho\left(\begin{bmatrix}
    A & BC_{\mK} \\
    B_{\mK}C & A_{\mK}
    \end{bmatrix}\right) < 1\right\}.
\end{equation}
Then, any such $\mK \in \mathcal{C}_q$ determines a dynamic policy $\pi_{\mK}$ with $\pi_{\mK}(\mathcal{H}_t) := \sum_{i=1}^t C_{\mK} A_{\mK}^{i - 1}B_{\mK} y_{t - i }$, where we set $\xi_0 = 0$.
Following the parameterization in (\ref{eq:closed-loop-param-policy-prob}), given the system plant dimension $n$, the policy optimization for LQG control becomes
\begin{equation}
\begin{aligned} \label{eq:LQG_policy-optimization}
    \min_{\mK} \quad & J_{\texttt{LQG}}(\mK) \\
    \mathrm{subject~to}\quad &  \mK \in \mathcal{C}_n. 
\end{aligned} 
\end{equation}
Throughout this survey, we make the standard assumption that $(A,B)$ is stabilizable and $(C,A)$ is detectable for the LTI system \eqref{eq:dynamics_discrete}, so that $\mathcal{C}_n$ is nonempty.

\subsubsection{$\mathcal{H}_\infty$-robust Control under Static State-feedback Policies}\label{subsubsect:state-feedback-robust-control}

The $\mathcal{H}_\infty$-norm of a closed-loop transfer function characterizes the worst case performance against adversarial disturbances $w_t$ with bounded energy. Considering the same dynamics in \eqref{eq:state-evolution}, then \eqref{eq:closed-loop-policy-prob} reduces to the $\mathcal{H}_\infty$ robust control problem: 
\begin{equation*} %
    \begin{aligned}
        \min_{\pi \in \Pi} \quad & J_\infty(\mathbf{u}):= \sup_{\|\mathbf{w}\|_{l_2} \leq 1} \lim_{T \to \infty} J_T(\mathbf{u}, \mathbf{w}, 0) \\
        \text{subject to} \quad & ~\eqref{eq:state-evolution}, \;
        u_t = \pi(\mathcal{H}_t),
    \end{aligned}
\end{equation*}
where we have assumed $x_0=0$ for simplicity.
Following (\ref{eq:closed-loop-param-policy-prob}) and by considering the same parameterization as LQR, we can equivalently express $J_\infty$ instead as a function on the set of static stabilizing policies $\mathcal{S}$, and thus the $\mathcal{H}_\infty$-robust policy problem becomes,
\begin{equation}
\begin{aligned} \label{eq:Hinf-policy-optimization}
    \min_{K} \quad & J_{\infty}(K) \\
    \mathrm{subject~to}\quad &  K \in \mathcal{S}. 
\end{aligned} 
\end{equation}
Note that \eqref{eq:Hinf-policy-optimization} is the state-feedback $\mathcal{H}_\infty$ control based on dynamics in \eqref{eq:state-evolution}, where the \textit{static linear} policy parameterization results in no loss of optimality \citep{zhou1988algebraic}. For the partially observed LTI system \eqref{eq:dynamics_discrete} where the state $x_t$ is not directly measured, we can consider the general output-feedback $\mathcal{H}_\infty$ control, in which a dynamic policy similar to \eqref{eq:controller_discrete} is required \citep{zhou1996robust}.

\section{Topology and Geometry of Stabilizing Policies}
\label{section:stablizing-policies}
Given an LTI system \eqref{eq:state-evolution}, represented with $(A,B) \in \Amatrices \times \Bmatrices$, the set of static stabilizing state-feedback policies $\mathcal{S} \subseteq \mathbb{R}^{m \times n}$ in  \eqref{eq:static-stabilizing-policies} has rich topological properties. For example, provided $(A,B)$ is stabilizable, we have $\mathcal{S} \neq \emptyset$ and any element of $\mathcal{S}$ can be identified with pole placement \citep{aastrom2021feedback}.   
Using the Jury stability criterion~\citep{jury1964theory} (the discrete-time version of the Routh–Hurwitz stability criterion \citep{parks1962new}), %
we can express $\mathcal{S}$ as a finite system of polynomial inequalities written in terms of the elements of $K \in \mathbb{R}^{m \times n}$. However, this set~of~polynomial inequalities is complicated and may not directly offer insights on topological properties of $\mathcal{S}$.

It is known that $\stableK$ is unbounded when $m \geq 2$, and its topological boundary $\partial \stableK = \{K \in \Kmatrices \;|\; \rho(A+BK) =1\}$ is a subset of $\Kmatrices$. Furthermore, as a result of the continuity of eigenvalues in the entries of the matrix, we can argue that $\stableK$ is an open set in $\Kmatrices$.
It can furthermore be shown that $\mathcal{S}$ is contractible. In particular, by the Lyapunov stability linear matrix inequality criterion and a Schur complement argument we can argue that:
\begin{fact} \label{fact:path-connected}
    The set of stabilizing static state-feedback policies $\mathcal{S}$ is path-connected.
\end{fact}
\noindent This property is vital for devising algorithmic iterates that have to reach a minimizer from any initial point in $\mathcal{S}$.

\subsection{Riemannian Geometry of Stabilizing Policies}\label{subsect:rie-geo}
Before diving into the geometry of stabilizing policies, we will first introduce basic concepts from Riemannian geometry. More specialized topics will be introduced in their respective sections in this paper. The starting point of departure 
for us is the realization that 
the set of stabilizing policies $\stableK$, as an open subset of $\mathbb{R}^{m \times n}$, is a smooth manifold; a geometric object, generically denoted by $\mathcal{M}$, that loosely speaking is locally Euclidean.
We call any smooth function $c:\mathbb{R} \to \mathcal{M}$ a smooth curve on $\mathcal{M}$. A tangent vector at $x \in \mathcal{M}$ is any vector $\dot{c}(0) = \frac{d}{dt}{c}(t)|_{t=0}$, where $c(\cdot)$ is a smooth curve passing through $x$ at $t=0$. The set of all such tangent vectors at $x$ is a vector space called the tangent space at $x$ denoted by $T_x \mathcal{M}$. Its dimension $\dim(T_x \mathcal{M})$ coincides with the dimension of the manifold. For open sets, such as $\stableK$, its tangent spaces identifies with the vector space it lies within: $T_K \stableK \equiv \mathbb{R}^{m \times n}$; this is referred to as the usual identification of the tangent space. The tangent bundle is simply the disjoint union of all tangent spaces $T \mathcal{M} := \{(x,v):x \in \mathcal{M}, v \in T_x\mathcal{M}\}$.

Let $F: \mathcal{M} \to \mathcal{N} \subset \mathbb{R}^M$ be a smooth function between two smooth manifolds. If we perturb $x \in \mathcal{M}$ along a direction $v \in T_x\mathcal{M}$, then the perturbation of the output from $F(\cdot)$ is another tangent vector in $T_{F(x)} \mathcal{N}$; the linear mapping $\diff F_x:T_x \mathcal{M} \to T_{F(x)} \mathcal{N}$ is called the differential of $F$ at $x$ and acts on any vector $v \in T_x \mathcal{M}$ as 
\begin{equation*}
    \diff F_x(v) := \frac{d}{dt} \big|_{t=0} (F \circ c)(t)
\end{equation*} 
where $c(\cdot)$ is any smooth curve passing $x$ at $t=0$ with velocity $v$. 

Considering the open set, and hence smooth manifold, of Schur stable matrices $\mathcal{A}$, we define the Lyapunov mapping $\lyap: \mathcal{A} \times \mathbb{R}^{n \times n} \to \mathbb{R}^{n \times n}$ that sends the pair $(A,Q)$ to the unique solution $P$ of 
\begin{equation}\label{eq:lyap-gen}
P = A P A^\transpose + Q.
\end{equation}
The following lemma is instrumental in geometric analysis of policies for linear systems.

\begin{lemma}\label[lemma]{lem:dlyap}
The differential of $\lyap$ at $(A,Q) \in \mathcal{A} \times \mathbb{R}^{n\times n}$ along $(E,F) \in T_{(A,Q)} (\mathcal{A} \times \mathbb{R}^{n\times n}) \equiv \mathbb{R}^{n \times n} \times \mathbb{R}^{n\times n}$ is
\begin{gather*}
    \diff \lyap_{(A,Q)}[E,F] =
    \lyap \big(A, E \lyap(A,Q) A^\transpose + A \lyap(A,Q) E^\transpose + F \big).
\end{gather*}
 For any $A \in \mathcal{A}$ and $Q, \Sigma \in \Amatrices$ we further have the so-called \emph{\lyaptrace} property,
 \[\tr{\lyap(A^\transpose,Q) \Sigma} = \tr{\lyap(A, \Sigma) Q}.\]
\end{lemma}

Moving on, a Riemannian metric $\langle \cdot,\cdot \rangle_x: T_x \mathcal{M} \times T_x \mathcal{M} \to \mathbb{R}$ on a smooth manifold $\mathcal{M}$ is an inner product that smoothly varies in $x$ with $x \in \mathcal{M}$. 
We call $(\mathcal{M},\langle \cdot,\cdot \rangle)$ a Riemannian manifold and often add a superscript to clarify specific Riemannian metrics. 
A (locally-defined) retraction is a smooth mapping $\mathcal{R}: \mathcal{T} \subset T\mathcal{M} \to \mathcal{M}$ where $\mathcal{T}$ contains an open subset of $(x,0_x) \in \mathcal{T}$ for each $x \in \mathcal{M}$, and the curve $c(t) := \mathcal{R}(x, t v)$ satisfies $c(0)=x$ and $\dot{c}(0)=v$.

The upshot of introducing
Riemannian geometry for policy optimization is now as follows. 
For any static feedback policy $K \in \stableK$, define the following Riemannian metric that depends on a solution of a Lyapunov equation\footnote{Compare with the Frobenius inner product $\tensor{V}{W}{}^\texttt{F} = \tr{V^\transpose W}$ which induces the so-called Euclidean geometry.}:
\begin{equation*}
    \tensor{V}{W}{K}^{\texttt{L}} \coloneqq \tr{V^\transpose  W  Y_K},
\end{equation*} where $Y_K := \lyap(A+ BK , \Sigma)$. In fact, this dependence varies smoothly in policy $K$ and thus we can show that $\tensor{\cdot}{\cdot}{}^\texttt{L}$ is a Riemannian metric on $\stableK$ referred to as the ``Lyapunov metric.''

\begin{theorem}\label{thm:Riemannian-metric}
If $\Sigma \succ 0$ then $(\stableK,\tensor{\cdot}{\cdot}{}^\textup{\texttt{L}})$ is a Riemannian manifold.\
\end{theorem}

\begin{figure}
    \centering
    \includegraphics[width = 0.45\textwidth]{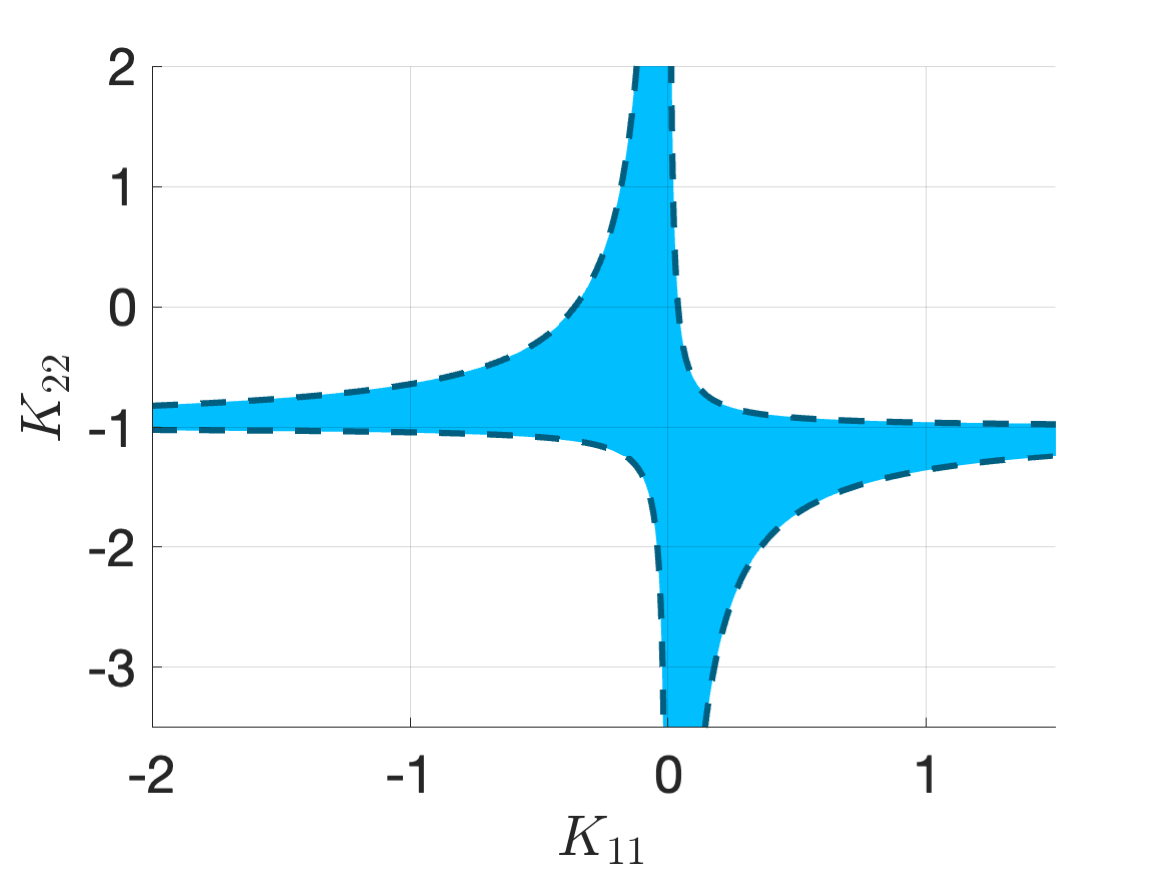}
    \caption{The 2-dimensional set of stabilizing policies subject to the off-diagonal sparsity constraint. The LTI system is $A = \begin{bmatrix}0.8 &1\\0 & 0.8\end{bmatrix}$ and $B = \begin{bmatrix}0 &1\\1 & 0\end{bmatrix}$. 
    }
    \label{fig:SLQR-landscape}
\end{figure}
The Riemannian machinery introduced for
stabilizing policies also allows addressing feedback design for certain classes of {\em constrained} policies. For example, when $\constraint \subset \mathbb{R}^{m \times n}$ and we restrict the policies to $\Theta = \substableK \coloneqq \stableK \cap \constraint$. From a \textit{topological} perspective, $\substableK \subset \stableK$ remains relatively open as $\stableK$ was open in $\Kmatrices$. However, in general, $\substableK$ is not only non-convex \citep{ackermann_parameter_1980} but also disconnected~\citep{feng_exponential_2019}. 
Nonetheless, one can show that sparsity constraint (which is of primary interest in network systems) and static output feedback policies lead to properly embedded submanifolds of $\stableK$~\citep{talebi_riemannian_2022}, thus entailing this Riemannian geometry as summarized in the following. See \Cref{fig:SLQR-landscape} for a numerical example of $\substableK$ with a sparsity constraint.

\begin{fact}
    For any sparsity constraint $\constraint_D \coloneqq \big\{K \in \Kmatrices \;|\; K_{i,j} = 0, \; (i,j) \not\in D\big\}$ with a index set $D \subset [m] \times [n]$, the set of sparse stabilizing policies $\substableK = \stableK \cap \constraint_D \subset \stableK$ is a properly embedded submanifold of dimension $|D|$. Also, each tangent space at $K \in \substableK$ identifies with $T_K\substableK \cong \constraint_D$.
\end{fact}

\begin{fact}
    For any output feedback constraint $\constraint_C \coloneqq \big\{K \in \Kmatrices \;|\; K = L C, \; L \in \Lmatrices \big\}$ with a full-rank output matrix $C$, the set of static output-feedback policies $\substableK = \stableK \cap \constraint_C$ is a properly embedded submanifold of $\stableK$ with dimension $m d$. Also, each tangent space at $K \in \substableK$ identifies with $T_K\substableK \cong \constraint_C$.
\end{fact}
\noindent We later see the implications of these fact when we study the (Riemannian) gradient and Hessian of a smooth cost in \S\ref{sec:lin-constrained-policy-performance}.

\subsection{Topology and Geometry of Dynamic Output-feedback Polices}

Herein, we focus on the output-feedback problem setup introduced in \S\ref{subsubsect:LQG}. Our parameterized family of policies will be the set of full-order dynamic output feedback policies $\Theta := \mathcal{C}_n$ in (\ref{eq:internallystabilizing}), and we next discuss some of its topological and geometrical properties.

\begin{figure}[t]
\centering
    \includegraphics[width = 0.4\textwidth]{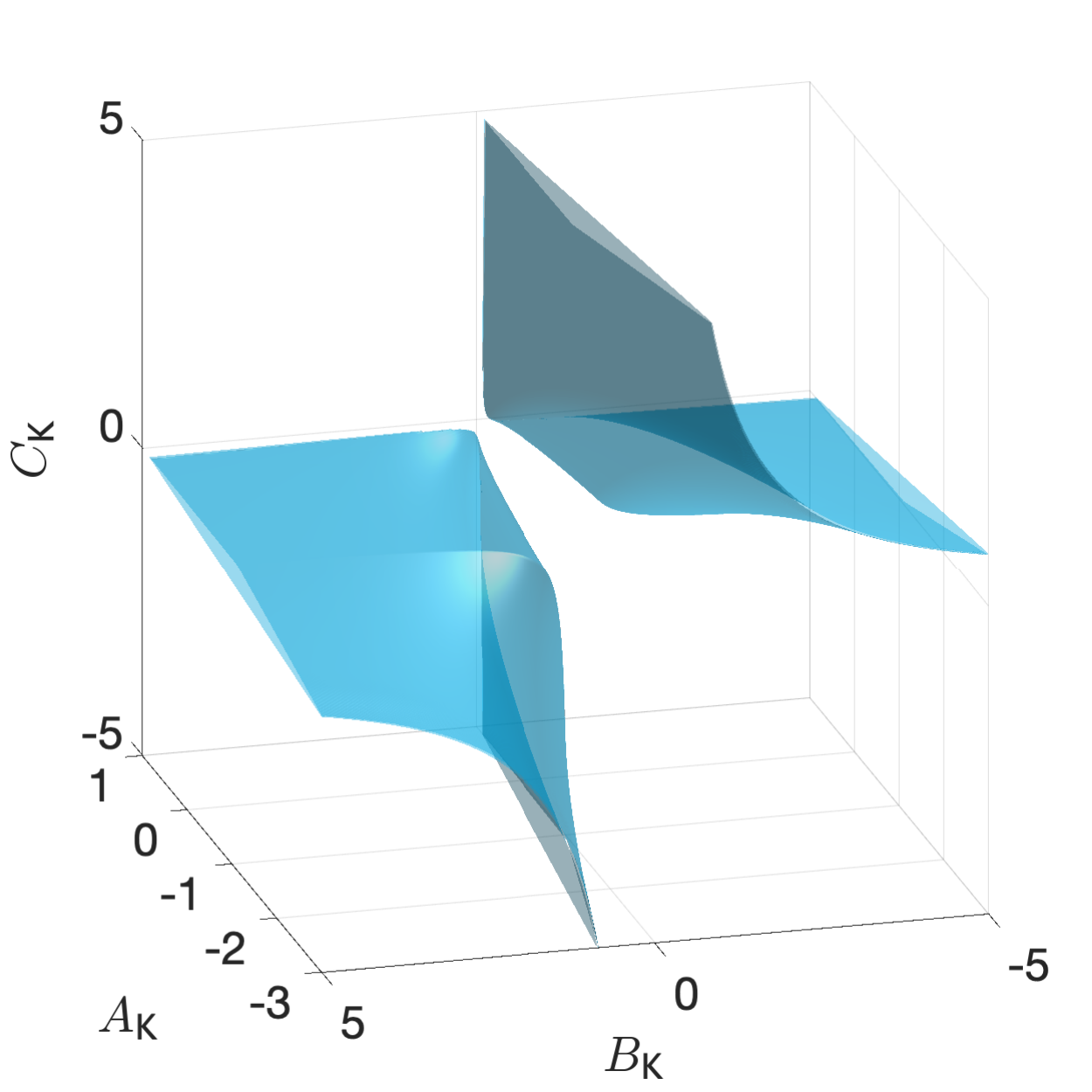}
    \hspace{15mm}
    \includegraphics[width=0.4\textwidth]{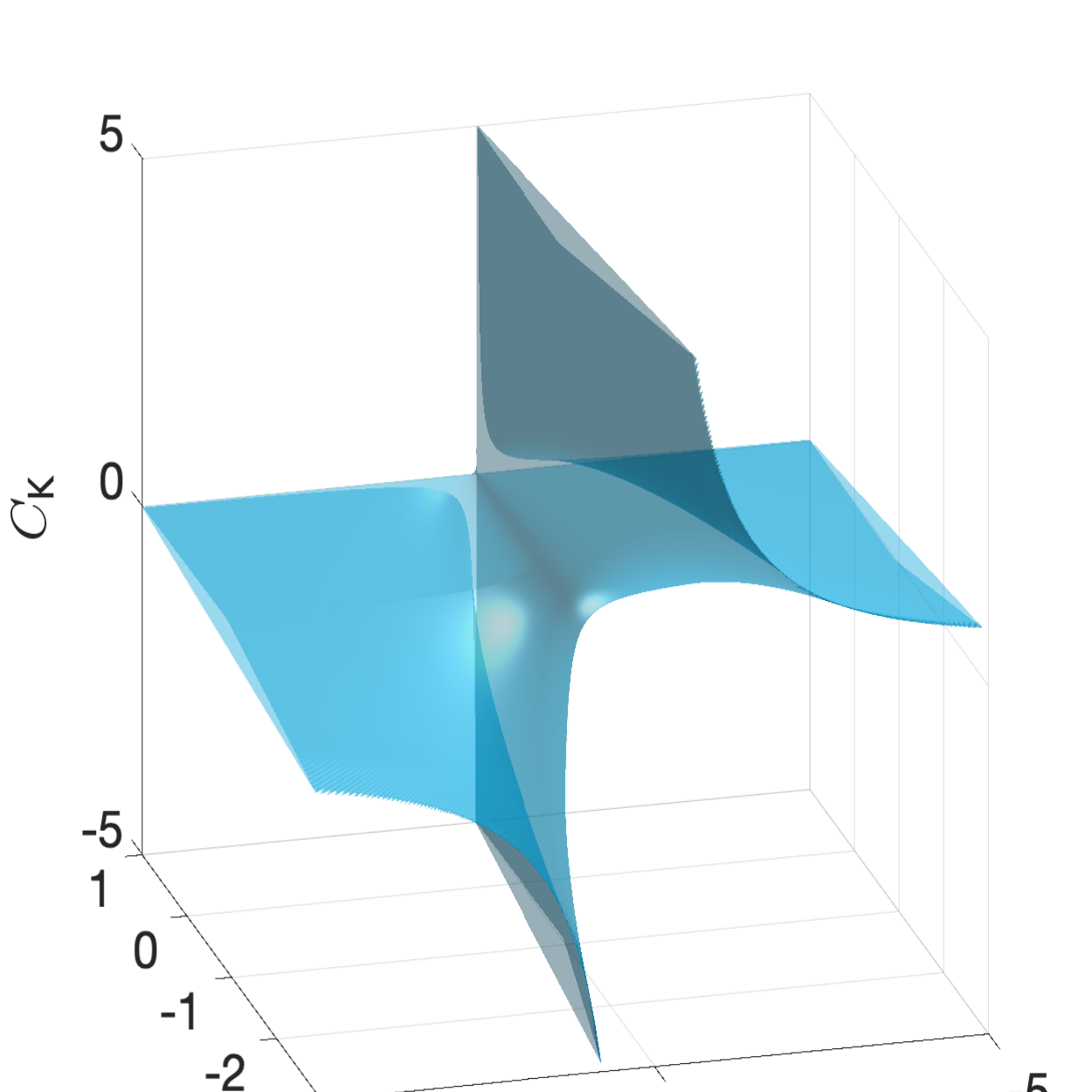}

    \caption{Illustration of the set of dynamic stabilizing policies $\mathcal{C}_1$ for an LTI system with $B = C = 1$ and: (a) with $A=1.1$ resulting in two path-connected components; (b) with $A=0.9$ resulting in a single path-connected component.}
    \label{fig:feasible_region}
\end{figure}

Similar to the static case, $\mathcal{C}_n$ is a nonconvex set. This is illustrated in \Cref{fig:feasible_region} with a numerical example.
It is also known that the set $\mathcal{C}_n$ is open and unbounded. 
In addition to the non-convexity, the set $\mathcal{C}_n$ can even be disconnected but has at most two diffeomorphic components that are captured by the following notion of ``similarity transformations in control'': for any invertible matrix $T \in \mathrm{GL}_n$, define the mapping $\mathscr{T}_T$ that sends any dynamic policy $\mK$ to
\begin{equation}\label{eq:def_sim_transform}
\mathscr{T}_T(\mK)
\coloneqq 
\begin{bmatrix}
0 & C_{\mK}T^{-1} \\
TB_{\mK} & TA_{\mK}T^{-1}
\end{bmatrix}.
\end{equation}
Note that the policy $\mK\in\mathcal{C}_n$ if and only if the transformed policy $\mathscr{T}_T(\mK) \in \mathcal{C}_n$. Indeed, the map
$
 \mK\mapsto \mathscr{T}_T(\mK)
$
is a diffeomorphism from $\mathcal{C}_n$ to itself for any such invertible matrix $T$.  \footnote{Note that the same input-output behavior of a policy \eqref{eq:controller_discrete} can be represented using different state-space models, e.g., using different coordinates for the internal policy state which is precisely captured by this similarity transformation.}

\begin{fact} \label{Fact:disconnectivity}
The set ${\mathcal{C}}_{n}$ has at most two path-connected components.  
{If ${\mathcal{C}}_{n}$ has two path-connected components $\mathcal{C}_n^{(1)}$ and $\mathcal{C}_n^{(2)}$, then $\mathcal{C}_n^{(1)}$ and $\mathcal{C}_n^{(2)}$ are diffeomorphic under the mapping $\mathscr{T}_T$, for any invertible matrix $T\in\mathbb{R}^{n\times n}$ with $\det T<0$. 
}
\end{fact}
The potential disconnectivity of ${\mathcal{C}}_{n}$ comes from the fact that the set of real invertible matrices $\mathrm{GL}_n=\{\Pi\in\mathbb{R}^{n\times n}\mid\,\det \Pi\neq 0\}$ has two path-connected components:
$
\mathrm{GL}^+_n=\{\Pi\in\mathbb{R}^{n\times n}\mid\,\det \Pi> 0\},
\,
\mathrm{GL}^-_n=\{\Pi\in\mathbb{R}^{n\times n}\mid\,\det \Pi< 0\}.
$ In other words, the nature of similarity transformations embedded in dynamic feedback policies may cause $\mathcal{C}_n$ to be disconnected.
The following results provide conditions that ensure $\mathcal{C}_n$ to be a single path-connected component.  

\begin{fact} %
If there exists a reduced-order stabilizing policy for \eqref{eq:dynamics_discrete}, \emph{i.e.}, ${\mathcal{C}}_{n-1} \neq \varnothing$, then $\mathcal{C}_n$ is path-connected. The converse also holds for systems with single-input or single-output, i.e., when $m=1$ or $p=1$ in \eqref{eq:dynamics_discrete}.
\end{fact}

We can immediately deduce the following facts: 1) For any open-loop unstable first-order dynamical system, \emph{i.e.}, $n = 1$ and $A > 0$, there exist no reduced-order stabilizing policies, \emph{i.e.},   $\mathcal{C}_{n-1} = \varnothing$; thus its associated set of stabilizing policies $\mathcal{C}_n$ must be disconnected. 2) For any open-loop stable systems, i.e., when $A$ is stable, we naturally have a reduced-order stabilizing policy, and thus the corresponding set of stabilizing policies is always path-connected. 
\Cref{fig:feasible_region} provides numerical examples for each case.

\subsection{Symmetries of Dynamic Output-feedback Policies: a Quotient Geometry}\label{subsubsec:invariants_of_feedback_policies}

\begin{figure}
    \centering
    \includegraphics[width=0.45\textwidth]{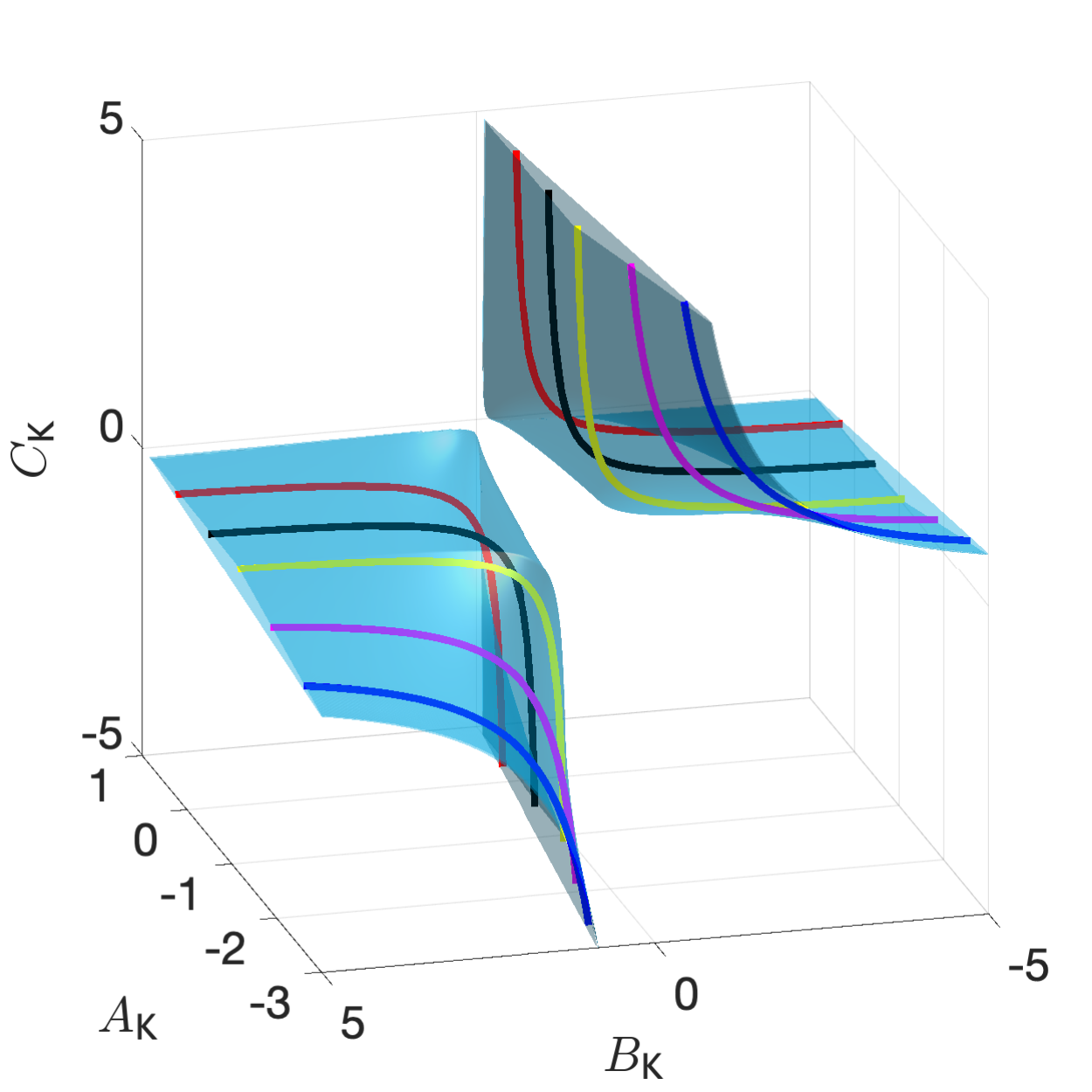}
    \caption{The region of stabilizing dynamic feedback policies $\mathcal{C}_1 \subset \mathbb{R}^3$ for the plant $(A,B,C)=(1.1,1,1)$. Each colored curve (red, purple, yellow, magenta, dark blue) is an individual orbit of policies. Note that $\mathcal{C}_1$ has 2 path-connected components.}
    \label{fig:controller_orbit_space}
\end{figure}

An alternative parameterization for dynamic output-feedback systems (\ref{eq:internallystabilizing}) is through turning $\mathcal{C}_n$ into a Riemannian quotient manifold of lower dimension. 
To see this, note that the group of similarity transformations $\{\mathscr{T}_T(\cdot):T \in \mathrm{GL}_n\} \equiv \mathrm{GL}_n$ is a group action acting smoothly on $\mathcal{C}_n$. Recall, we are treating the open subset $\mathcal{C}_n \subset \mathbb{R}^{(n + m) \times (n + p)}$ as a smooth manifold. 
The orbit of $\mK \in \mathcal{C}_n$ is then the collection of controllers reachable from $\mK$ via similarity transformation:
\begin{equation*}
    [\mK]:=\left\{\begin{bmatrix} 0 & C_{\mK}T^{-1} \\ TB_{\mK} & TA_{\mK}T^{-1} \end{bmatrix}: T \in \mathrm{GL}_n \right\}.
\end{equation*}
A few examples of such orbits are shown in Figure \ref{fig:controller_orbit_space} in distinct colors. Recall that LQG cost is constant on each orbit; this serves as the basis for the PO for LQG over the so-called quotient space.

The quotient set (also known as the orbit set) of $\mathcal{C}_n$ is simply the collection of all orbits:
\begin{equation*}
    \mathcal{C}_n/\mathrm{GL}_n := \{[\mK] : \mK \in \mathcal{C}_n\}.
\end{equation*} 
We equip $\mathcal{C}_n/\mathrm{GL}_n$ with the induced quotient topology, defined as the finest topology in which the quotient map $\pi:\mathcal{C}_n \to \mathcal{C}_n/\mathrm{GL}_n,$ sending each $\mK$ to the orbit $[\mK]$, is continuous. The resulting topological space is called the quotient space.

The next step is to design a smooth structure for $\mathcal{C}_n/\mathrm{GL}_n$, turning it into a smooth manifold. For arbitrary quotient spaces, if there exists a smooth structure in which the quotient map is a smooth submersion, we call the resulting quotient space a smooth quotient manifold. In this context, the original smooth manifold is called the total manifold. %
Unfortunately, quotient spaces are often not even Hausdorff. Recall a  topological space is Hausdorff when any pair of points can be separated into disjoint neighborhoods of the corresponding points; this property implies all sequences have a unique limit. Therefore for non-Hausdorff quotient spaces, optimization is hopeless because limits cannot even be defined!

As it turns out, $\mathcal{C}_n/\mathrm{GL}_n$ is non-Hausdorff. The reason is the existence of non-controllable and non-observable, yet stabilizing, dynamic controllers which acts as ``jumps.'' This is explained in more detail in \cite{hazewinkel_moduli_1976}. Fortunately, the quotient space $\mathcal{C}_n^{\min}/\mathrm{GL}_n$ is Hausdorff \cite[Lemma 4.1]{kraisler2024output},
where $\mathcal{C}_n^{\min}$ are all full-order \textit{minimal} policies--- policies with controllable and observable state space form. This follows from the remarkable theorem in \cite{hazewinkel_moduli_1976} proving the orbit space of minimal linear systems admits a smooth quotient manifold structure. It can be shown $\dim(\mathcal{C}_n^{\min}/\mathrm{GL}_n)=nm+np$, an order of magnitude smaller. So, in the context of smooth optimization, we have a significantly smaller search space.

Before we continue, let us discuss the tangent space of $\mathcal{C}_n$. As an open subset, the tangent space of $\mathcal{C}_n$ coincides with its linear span: tangent \textit{vectors} to the open set $\mathcal{C}_n$ are simply 
\begin{equation*}
    \mathbf{V} = \begin{bmatrix} 0 & G \\ F & E\end{bmatrix}
\end{equation*} for any matrices $E \in \mathbb{R}^{n \times n}$, $F \in \mathbb{R}^{n \times p}$, and $G \in \mathbb{R}^{m \times n}$. The resulting vector space will be denoted $\mathcal{V}_n$; so, we write $T_K \mathcal{C}_n \equiv \mathcal{V}_n$.

In order to adopt smooth optimization techniques for $\mathcal{C}_n^{\min}/\mathrm{GL}_n$, we at last must equip the smooth quotient manifold with a retraction and Riemannian metric. To do this, we will discuss a correspondence of these two constructs between the total manifold $\mathcal{C}_n^{\min}$ and the quotient manifold $\mathcal{C}_n^{\min}/\mathrm{GL}_n$. We can show that there is an invertible correspondence between Riemannian metrics on $\mathcal{C}_n^{\min}/\mathrm{GL}_n$ and \textit{similarity-invariant} Riemannian metrics on $\mathcal{C}_n^{\min}$, by which we mean
\begin{equation}
    \langle \mathbf{V}, \mathbf{W} \rangle_{\mK} = \langle \mathscr{T}_T(\mathbf{V}), \mathscr{T}_T(\mathbf{W}) \rangle_{\mathscr{T}_T(\mK)} \label{eq:metric-inv}.
\end{equation} A related correspondence holds for retractions. See \cite[\S9]{boumal_introduction_2023} and \cite[\S4]{kraisler2024output} for how to induce such a Riemannian metric and retraction onto the quotient manifold.

Now we will define a similarity-invariant Riemannian metric satisfying (\ref{eq:metric-inv}). Let $A_\cl(\mK)$,$B_\cl(\mK)$, and $C_\cl(\mK)$ denote the matrices corresponding to the (augmented) closed-loop system in (\ref{eq:closed_loop_system}). A consequence of the minimality of any $\mK \in \mathcal{C}_n^{\min}$ is (1) $A_\cl(\mK)$ is Hurwitz and (2) the closed loop system $(A_\cl(\mK),B_\cl(\mK),C_\cl(\mK))$ is also minimal. The latter follows from the Popov-Belevitch-Hautus test. As a result, the controllability and observability Grammains of the closed-loop system satisfy 
\begin{subequations}
    \begin{align}\label{eq:pos_def_grammians}
        \mathcal{W}_c(\mK) &:= \lyap(A_\cl(\mK), B_\cl(\mK)B_\cl(\mK)^\transpose ) > 0\\
        \mathcal{W}_o(\mK) &:= \lyap(A_\cl(\mK), C_\cl(\mK)^\transpose C_\cl(\mK)) > 0.
    \end{align}
\end{subequations} 
Now, we introduced the so-called ``Krishnaprasad-Martin (KM) metric'' defined as
\begin{subequations}\label{eq:KM-metric}
\begin{align}
    \langle \mathbf{V}_1, \mathbf{V}_2 \rangle_\mK^{\texttt{KM}} &:= w_1 \text{tr}[\mathcal{W}_o(\mK) \cdot \mathbf{E}(\mathbf{V}_1) \cdot \mathcal{W}_c(K) \cdot \mathbf{E}(\mathbf{V}_2)^T] \\
    &\qquad + w_2 \text{tr}[\mathbf{F}(\mathbf{V}_1)^T \cdot \mathcal{W}_o(\mK) \cdot \mathbf{F}(\mathbf{V}_2)] \\
    &\qquad + w_3 \text{tr}[\mathbf{G}(\mathbf{V}_1) \cdot \mathcal{W}_c(\mK) \cdot \mathbf{F}(\mathbf{V}_2)^T] 
\end{align}
\end{subequations} where $w_1>0$ and $w_2,w_3 \geq 0$ are design constants and 
\begin{align*}
    \mathbf{E}(\mathbf{V}) := \begin{bmatrix} 0 & BG \\ FC & E \end{bmatrix}, \quad
    \mathbf{F}(\mathbf{V}) := \begin{bmatrix} 0 & 0 \\ 0 & F\end{bmatrix}, \quad 
    \mathbf{G}(\mathbf{V}) := \begin{bmatrix} 0 & 0 \\ 0 & G \end{bmatrix}.
\end{align*}
We can show that the inner-product in \eqref{eq:KM-metric} varies smoothly in $\mK$ and, by referring to (\ref{eq:pos_def_grammians}) for each $K\in\mathcal{C}_n^{\min}$, establish the following result (analogous to the Riemannian metric in \Cref{thm:Riemannian-metric}).

\begin{theorem}\label{thm:Riem-metric-dynamic}
   The mapping that sends each $\mK \in \mathcal{C}_n^{\min}$ to the inner-product $\langle ., . \rangle^{\texttt{KM}}_\mK$ induces a Riemannian metric on $\mathcal{C}_n^{\min}$ which is similarity-invariant; i.e., satisfies (\ref{eq:metric-inv}).
\end{theorem}

\section{Geometry of Performance Objectives on Stabilizing Policies} \label{section:LQ-performance}

In this section, we turn our attention toward performances measures described in \S\ref{parameterization} and the interplay of domain geometry and the landscape of the cost functions associated with each performance measure. Before getting to each specific metric, we review some of the geometric constructs that are essential in characterizing the first and second-order variations of smooth cost functions, subsequently used for optimization.

\subsection{Riemannian Geometry and Policy Optimization}

For brevity, we introduce these constructs for the smooth manifold $\stableK$ which also holds for any other smooth manifold. A vector field $V:\stableK \to T\stableK$ is a mapping that smoothly assigns every $K \in \stableK$ to a tangent vector $V_K \in T_K \stableK$. A vector field induces a mapping on the space of smooth functions, sending $J:\stableK \to \mathbb{R}$ to $VJ:\stableK \to \mathbb{R}$ as follows: 
\begin{equation}
    VJ(K) := \diff J_K(V_K)
\end{equation} 

In this section, to distinguish tangent vectors from vector fields, we will denote the former with $V_K$ to emphasize that $V_K$ is a tangent vector in $T_K \stableK$. Now, we can define the (Riemannian) gradient of $J$ with respect to the Riemannian metric $\tensor{\cdot}{\cdot}{}^\texttt{L}$ on $\stableK$, denoted by $\grad{J}$. In particular, $\grad J$ is the unique vector field satisfying
\begin{equation}\label{eq:grad-def}
    \tensor{V}{\grad J}{}^\texttt{L} = V J,
\end{equation} for all vector fields $V$.

In order to define second-order variations of a smooth functions $J$, such as the Riemannian Hessian, we must introduce a notion of directional derivatives on manifolds. This is referred to as the \textit{(affine) connection}, denoted by $\nabla$. Consider two vector fields $V$ and $W$ on $\stableK$. Then, the connection $\nabla$ allows us to define $\nabla_V W$, which itself is a vector field, at each $K \in \stableK$ as the directional derivative of $W$ along $V_K \in T_K \stableK$. Each Connection on $\stableK$ is uniquely identified with $\dim(\stableK)^3$ number of smooth functions on $\stableK$ called the \textit{Christoffel Symbols}. With the Christoffel Symbols associated with $\nabla$, in order to compute $(\nabla_V W)_K$, we only need the direction $V_K$ and the vector field $W$ locally; we do not need other evaluations of $V$. 

Every Riemannian metric uniquely induces a \textit{compatible} affine connection known as the \textit{Levi-Civita connection}. It is the unique connection which satisfies
\begin{equation}
    \nabla_U \tensor{V}{W}{} = \tensor{\nabla_U V}{W}{} + \tensor{V}{\nabla_U W}{}
\end{equation} for any vector fields $U,V,W$ of $\stableK$.\footnote{To be precise, for technical reasons regarding the uniqueness we must also require that $\nabla$ is ``symmetric;'' see \cite{talebi_policy_2023}. } Hereafter, we let $\nabla$ and $\overline{\nabla}$ denote the Levi-Civita connections of $\stableK$ compatible with $\langle .,. \rangle^\texttt{L}$ and the Frobenius inner product $\langle V, W\rangle^\texttt{F} := \tr{V^\transpose W}$, respectively. %
The macron in $\overline{\nabla}$ indicates the ``flatness'' of this Euclidean directional derivative.

Letting $\fX(\stableK)$ denote the family of vector fields on $\stableK$, we can define the Riemannian Hessian of $J$ in two equivalent ways:

\begin{subequations}\label{eq:hess-def}
\begin{align}
    \hess J(V) &:= \nabla_V \grad J \\
    \hess J(V, W) &:= \langle \nabla_V \grad J, W \rangle^\texttt{L}
\end{align}
\end{subequations}
for any $V,W \in \fX(\stableK)$.
Both forms can be used interchangeably. The former is used in the context of Riemannian-Newton optimization. The latter is used in the context of defining a matrix representing the Hessian matrix of $f$. It should be noted that for both definitions, $V$ and $W$ do not have to be tangent vector \textit{fields}; they could simply be tangent vectors. To emphasize our evaluation at a specific $K \in \stableK$, we will write $\hess J |_K $ for both definitions.

As usual, the Euclidean Hessian is equivalently defined as $\overline{\nabla}^2 J(V) := \overline{\nabla}_V \overline{\nabla} J$.
At last, we also introduce an atypical yet efficient notion known as the pseudo-Euclidean Hessian which essentially ignores the curvature of the manifold in quantifying the second order behavior of a smooth function. It is constructed by applying the \textit{Euclidean affine connection} $\overline{\nabla}$ on the (Riemannian) gradient $\grad J$ as follows:
\begin{equation}
    \overline{\hess} J(V) := \overline{\nabla}_V \grad J
\end{equation}
which will be compared
 with the Riemannian Hessian, denoted by $\hess J$.

\definecolor{myblue}{RGB}{0, 190, 255}
\definecolor{myviolet}{RGB}{126, 47, 142}
\begin{figure}[t]
    \centering
    \includegraphics[width = 0.35\textwidth, trim ={4cm  2cm 4cm 2cm},clip]{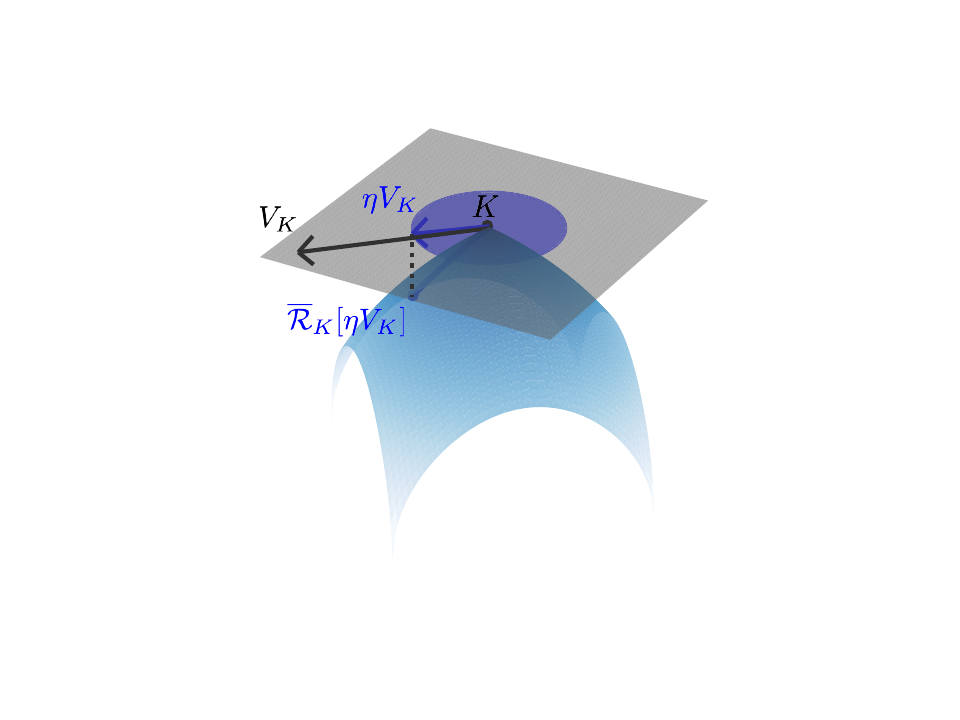}
    \caption{Local retraction defined by the stability certificate: A schematics of the gray plane exemplifying a tangent space at a point $K$ on the blue manifold. The stability certificate provides a (purple) neighborhood of the origin (in every tangent space) such that an efficient ``local retraction'' $\overline{\mathcal{R}}_K$can be obtained, such that every tangent vector $V_K$ can be ``retracted'' to $\overline{\mathcal{R}}_K[\eta V_K]$ after proper scaling by the stability certificate $s_K$ with $\eta := s_K(V_K)$. 
    }
    \label{fig:stabcert}
\end{figure}

\subsection{Stability Certificate for the Euclidean Retraction}
For the open submanifold of a vector space, such as the static feedback policies $\mathcal{S}$, a useful example of a retraction is the Euclidean retraction: $\overline{\mathcal{R}}_K(V_K):= K + V_K$, which is computationally efficient. We emphasize that this is not well-defined globally on $T_K\stableK$, and thus motivates us to further determine the local neighborhood on which it will be well-defined; see \Cref{fig:stabcert}.

\begin{lemma}\label[lemma]{lem:stability-cert}
For any direction $V_K \in T_K\stableK \cong \Kmatrices$ at any point $K \in \stableK$, if 
\begin{equation*}
    0 \leq \eta \leq s_K(V_K) \coloneqq \frac{1}{2 \lambda_{\max}\left(\lyap(A_{\cl}^\transpose, I)\right) \|BV_K\|_2 }
\end{equation*}
where $A_{\cl} = A+BK$, then $\overline{\mathcal{R}}_K(\eta V_K) = K + \eta V_K \in \stableK$.
\end{lemma}

A few remarks are in order. First, since $\stableK \subset \Kmatrices$ is open, the stability certificate offers a closed-form expression of a continuous lower bound on the radius to instability: $s_K(V_K) \leq \sup \{t:t > 0 , K + tV_K \in \stableK\}$, which has an unknown closed-form expression. In the structured LQR setup, since $\constraint$ is an affine space, $\widetilde{\stableK}$ is relatively open and thus $s_K(\cdot)$ is also a stability certificate on $\widetilde{\stableK}$. Given $K \in \widetilde{S}$ and $V_K\in T_K \widetilde{\stableK}$, then
\begin{equation}\label{eq:RGD-step}
    K^+ := \overline{\mathcal{R}}_K(\eta V_K) = K + \eta V_K
\end{equation} for $\eta = \min(1, s_K(V_K))$ renders $K^+ \in \widetilde{S}$. Thus, for iterative update of policy $K$ in (\ref{eq:RGD-step}), the chosen step size guarantees feasibility and stability of $K^+ \in \widetilde{S}$. 

For the LQR setup, the direction $V_K$ could be the negated Riemannian gradient $-\grad J_\texttt{LQR}(K)$ or Euclidean gradient $-\overline{\nabla} J_\texttt{LQR}(K)$ of the LQR cost. The direction could also incorporate second-order information, in the form of a Riemannian-Newton optimization, formulated as the solution $V_K \in T_K\stableK$ of any of these linear equations:
\begin{align*}
    \hess J_\texttt{LQR} |_K (V) &= -\grad J_\texttt{LQR}(K) \\
    \overline{\nabla}^2 J_\texttt{LQR} |_K(V) &= -  \overline{\nabla} J_\texttt{LQR}(K) \\
    \overline{\hess} J_\texttt{LQR} |_K(V) &= -  \grad J_\texttt{LQR}(K).
\end{align*} Additionally, for structured LQR setup, we show that these first and second order variations can be obtained similarly for the constrained cost $\widetilde{J} = J |_{\substableK}$ using \Cref{thm:submanifold} and the policy update proceeds similarly.

\begin{remark}\label{rem:hewer}
In the absence of constraint when $\substableK=\stableK$, the Hewer's algorithm introduces the following updates 
$K^+ = -(R + B^\transpose P_{K} B)^{-1}B^\transpose P_{K
} A$ with $P_K = \lyap(A_{\cl}^\transpose, Q+ K^\transpose R K)$,
which is shown to converge to the global optimum quadratically. Somewhat interestingly, it
in can be written as 
$K^+ = K + \widehat{V},$
with a ``Riemannian quasi-Newton'' direction $\widehat{V}_t$ satisfying,
\[\widehat{H}_{K}\widehat{V} = - \grad{J}_{\texttt{LQR}}(K),\]
where $\widehat{H}_{K}: = R + B^\transpose P_{K} B$ is a positive definite approximation of $\hess J_\texttt{LQR}|_K$ and $\overline{\hess} J_\texttt{LQR}|_K$.
The \emph{algebraic coincidence} is that the unit stepsize remains stabilizing throughout these quasi-Newton updates. 
We will also see that for the unconstrained LQR problem a small enough (fixed) step-size is sufficient.
In general though, we do not expect such step sizes to be stabilizing, and particularly on constrained submanifolds $\substableK$ one needs to instead utilize
the stability certificate developed in \Cref{lem:stability-cert} to guarantee the stability of each policy iterate.
\end{remark}

\subsection{Linear Quadratic Regulator (LQR)}

We here discuss the non-convex geometry in the policy optimization for LQR. We consider both the standard (unconstrained) case, where the policy parameter $K$ can be a dense matrix, and the constrained case, where $K$ has extra linear constraints in addition to the stabilizing requirement--such as sparsity or output measurement.

\subsubsection{Unconstrained Case}\label{susec:standard-LQR-perf}

First, by the \lyaptrace~ property in \Cref{lem:dlyap}, one can show that for the static feedback parameterization $u_t = K x_t$ where $K \in \stableK$, both the cost functions in \eqref{eq:LQR-noise} and \eqref{eq:LQR-noiseless} are equivalent with 
\begin{equation} \label{eq:LQR-cost-lyapunov}
J_{\texttt{LQR}}(K) = \frac{1}{2} \tr{P_K \Sigma} = \frac{1}{2} \tr{(Q+K^\transpose R K) Y_K},
\end{equation}
where $P_K = \lyap(A_{\cl}^\transpose, Q+ K^\transpose R K)$, $Y_K = \lyap(A_{\cl}, \Sigma)$, and $A_{\cl} = A+BK$.
Next, the first and second variations of the smooth cost $J_{\texttt{LQR}} \in C^\infty(\stableK)$ with respect to the Riemannian and Euclidean metrics are obtained in the following proposition.

\begin{proposition}
    Consider the Riemannian manifold $(\stableK, \tensor{.}{.}{}^\textup{\texttt{L}})$. Then $J_{\texttt{LQR}}(\cdot)$ is smooth and
    \begin{align*}
        \grad J_\texttt{LQR}(K) &= RK + B^\transpose P_K A_{\cl} \\
        \overline{\nabla} J_\texttt{LQR}(K) &= (RK + B^\transpose P_K A_{\cl}) Y_K.
    \end{align*} See \cite{talebi_policy_2023} for explicit formulae of the Riemannian Hessian, and other second order variations of $J_\texttt{LQR}$.

\end{proposition}

Next, we review the properties of this cost function critical to policy optimization.
\begin{lemma}\label[lemma]{lem:coercive}
Suppose $\Sigma, Q, R \succ 0$ and $(A,B)$ is stabilizable. Then, the function $J_{\texttt{LQR}}(\cdot)\colon\mathcal{S} \to \bR$,
\begin{itemize}
    \item[(a)] is real analytic, and in particular smooth.
    
    \item[(b)] is coercive: $K \to \partial \stableK$ or $\|K\|_F \to \infty$ implies $J_\texttt{LQR}(K) \to \infty$;

    \item[(c)] admits a unique global minimum $K^*$ on $\stableK$ satisfying $K^* = -(R + B^\transpose P_K B)^{-1} B^\transpose P_K A$;
    
    \item[(d)] has compact sublevel sets $\stableK_\alpha := \{K:J_\texttt{LQR}(K) \leq \alpha\}$ for each finite $\alpha$;

    \item[(e)] is gradient dominant on each sublevel set: there exists a constant $c_1 >0$ such that
    \[c_1 [J_{\texttt{LQR}}(K)-J_{\texttt{LQR}}(K^*)] \leq  \| \overline{\nabla} J_{\texttt{LQR}}(K) \|_F^2, \quad \forall K \in \stableK_\alpha;\]

    \item[(f)] has $L$-Lipchitz gradient on each sublevel set: there exists a constant $L>0$ such that
    \[\|\overline{\nabla} J_{\texttt{LQR}}(K) - \overline{\nabla} J_{\texttt{LQR}}(K')\|_F \leq L  \|K - K'\|_F, \quad \forall K,K' \in \stableK_\alpha;\]

    \item[(g)]  admits lower and upper quadratic models on each sublevel set: there exists constants $c_2 >0$ and $c_3 >0$ such that
    \[c_2 \|K-K^*\|^2_F \leq  J_{\texttt{LQR}}(K)-J_{\texttt{LQR}}(K^*) \leq c_3 \|K-K^*\|^2_F, \quad \forall K \in \stableK_\alpha.\]
\end{itemize}
\end{lemma}

These properties of $J_\texttt{LQR}$ are quite essential in providing theoretical guarantees for different optimization schemes. In particular, the smoothness, Lipschitz continuity, and quadratic models are common in convex optimization whereas the gradient dominance enables global convergence guarantees despite non-convexity of $J_{\texttt{LQR}}$ in $K$. Finally, note that these properties holds on each (fixed) sublevel set of the cost which often is chosen to contain the initial policy $K_0$. These has been made possible due to the coercive property of $J_{\texttt{LQR}}$ that results in compact sublevel sets.

\subsubsection{Linearly Constrained Case}\label{sec:lin-constrained-policy-performance}

Here, we will discuss the Riemannian geometry of the LQR cost in the structured LQR setup. Since $\widetilde{\stableK} := \stableK \cap \constraint \subset \stableK$ is an embedded submanifold, the Riemannian metric $\langle .,. \rangle^\texttt{L}$ can be equipped onto $\widetilde{\stableK}$ simply by restricting its domain $T\stableK$ onto $T\widetilde{\stableK}$. Also, for any smooth function $J$ on $\stableK$, we let $\widetilde{J}:=J|_{\widetilde{\stableK}}$ be its restriction to $\substableK$. However, the gradient and Hessian of $\widetilde{J}$ will not relate to those of $J$ so simply. As for gradient $\grad \widetilde{J}:\substableK \to T\substableK$, our Euclidean intuition is correct and so it can be related to $\grad J$ by the ``tangential projection'' operator $\pi^\top$---the generalization of orthogonal projection with respect to the Riemannian metric. On the other hand, for the Riemannian Hessian of this restricted cost, denoted by $\hess \widetilde{J}$, our Euclidean intuition fails as this correspondence cannot be explained by merely a projection operator. In fact, the curvature of the underlying Riemmanian manifold affects this second-order information. This can be captured precisely by the \textit{Weingarten map} $\mathbb{W}_U(V) \in \fX(\widetilde{\stableK})$ as the unique vector field satisfying
\begin{equation}
    \langle \mathbb{W}_U(V), W \rangle^\texttt{L} = \langle U, \nabla_V W - \pi^\top \nabla_V W \rangle^\texttt{L}
\end{equation} for all $W \in \fX(\widetilde{\stableK})$. These relations are summarized below; see \cite{talebi_policy_2023} for further details.

\begin{theorem}\label{thm:submanifold}
    Let $J:\stableK \to \mathbb{R}$ and $\widetilde{J}:=J|_{\widetilde{\stableK}}$. Then over $\widetilde{\stableK}$, we have
    \begin{equation}
        \grad \widetilde{J} = \pi^\top \grad J.
    \end{equation} Furthermore, for any $V \in \fX(\widetilde{\stableK})$, we have
    \begin{equation}
        \hess \widetilde{J}(V) = \pi^\top \hess J(V) + \mathbb{W}_U (V),
    \end{equation} where $U := \grad J - \pi^\top \grad J$. 
\end{theorem}

This result enables us to obtain explicit formulae for the Riemannian gradient and Hessian of any smooth costs on $\stableK$ when restricted to $\widetilde{\stableK} = \stableK \cap \constraint$. In fact, this is proved for any embedded submanifold of an abstract Riemannian manifold.  As expected, these geometric derivatives will be affected by the curvature of $(\stableK, \langle .,. \rangle^\texttt{L})$ which is accurately captured by the second term of the Weingarten mapping. In explicit form, these can be computed using the Christoffel symbols associated with induced Levi-Civita connection $\nabla$; see \cite[Prop. 3.4]{talebi_policy_2023} for the general proof and explicit formulae of these quantities.

By direct application of these results to the constraint LQR cost $\widetilde{J}_\texttt{LQR} := J_\texttt{LQR}|_{\widetilde{\stableK}}$, we can give explicit formulae for Riemannian gradient and Hessian:
\begin{align*}
    \grad \widetilde{J}_\texttt{LQR}(K) &= \pi^\top (RK + B^\transpose P_K A_{\cl}) \\
    \overline{\nabla} \widetilde{J}_\texttt{LQR}(K) &= \text{Proj}_K((RK + B^\transpose P_K A_{\cl})Y_K),
\end{align*} where $\text{Proj}_K$ is the ordinary Euclidean projection operator from $T_K \stableK$ onto $T_K \widetilde{\stableK}$.
Similarly, the Riemannian, Euclidean, and Pseudo-Euclidean Hessians for $\widetilde{J}_\texttt{LQR}(\cdot)$ can be obtained.

\begin{figure}
    \centering
    \subfigure[]{\includegraphics[width = 0.46\textwidth]{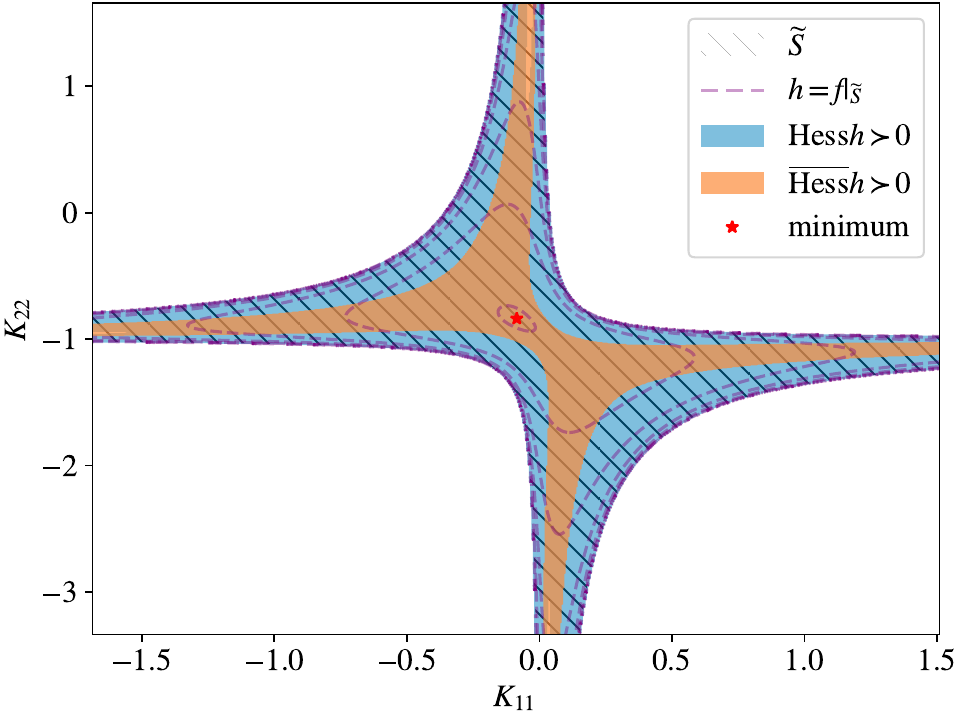}}
    \subfigure[]{ \includegraphics[width = 0.52\textwidth]{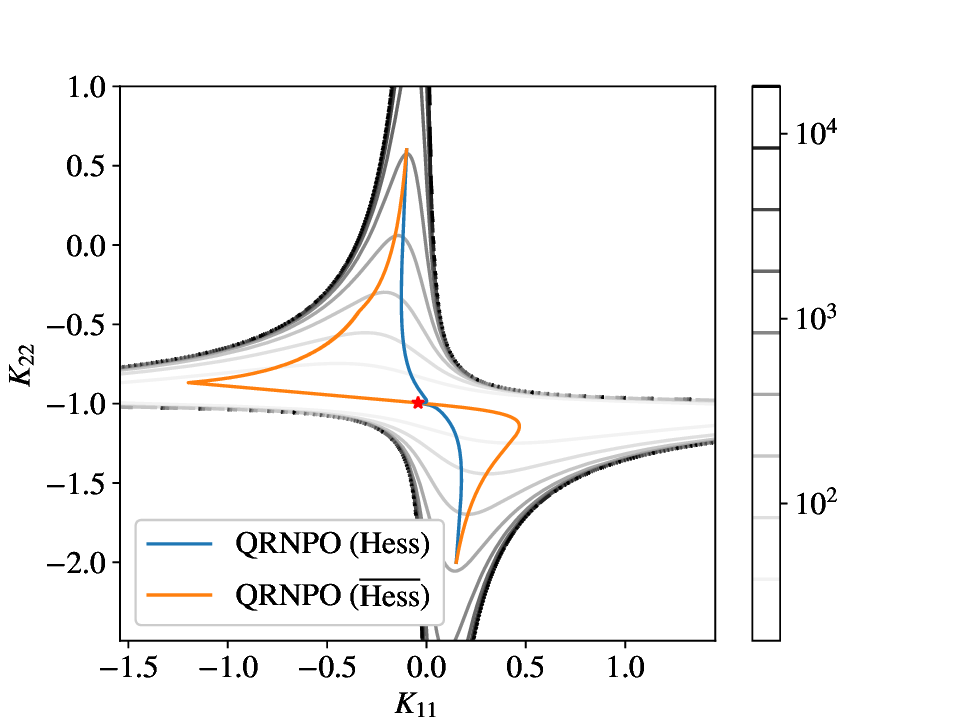}}
    \caption{A numerical example for the feasible set of diagonal stabilizing controller for the \ac{slqr} problem of \S 2. a) the region on which the Riemannian Hessian is positive definite versus that of Pseudo-Euclidean Hessian. (b) Trajectories of the QRNPO algorithm \citep{talebi_riemannian_2022} using Riemannian Hessian in blue, incorporating the curvature, versus the Euclidean Hessian in orange, where the former takes more efficient trajectories towards the global minimum (denoted by the red star).}
    \label{fig:LQR-landscape-constrained}
\end{figure}

In \Cref{fig:LQR-landscape-constrained}, we provide a numerical illustration of how this Riemannian metric is useful, and how the curvature information enables more
efficient algorithms.\footnote{The figure pertains to the example as in \Cref{fig:SLQR-landscape}, where the feedback gain is constrained to be diagonal.}

\subsection{Linear Quadratic Gaussian (LQG) Control}\label{subsec:LQG}

In this subsection, we move to discuss the geometry in policy optimization for LQG control \eqref{eq:LQG_policy-optimization}. As we will see below, while sharing certain similarities, the non-convex LQG landscape is richer and more complicated than the LQR case.

Recall that the set of stabilizing policies $\mathcal{C}_n$  has at most two path-connected components that are diffeomorphic to each other under a particular similarity transformation (Fact \ref{Fact:disconnectivity}). As similarity transformations do not change the input/output behavior of dynamic policies, it makes no difference to search over any path-connected component even if $\mathcal{C}_n$ is not path-connected. This feature brings positive news for local policy search algorithms over dynamic stabilizing policies $\mathcal{C}_n$. %

\subsubsection{Spurious Stationary Points and Global Optimality}

Similar to the LQR case, for any stabilizing policy $\mK\in\mathcal{C}_n$, the LQG cost function $J_{\texttt{LQG}}(\mK)$ in \eqref{eq:LQG_policy-optimization} has the following expressions: 
\begin{equation}\label{eq:LQG_cost_formulation_discrete}
J_{\texttt{LQG}}(\mK) = \operatorname{tr}
\left(
\begin{bmatrix}
Q & 0 \\ 0 & C_{\mK}^\transpose R C_{\mK}
\end{bmatrix} X_\mK\right)
=
\operatorname{tr}
\left(
\begin{bmatrix}
W & 0 \\ 0 & B_{\mK} V B_{\mK}^\transpose
\end{bmatrix} Y_\mK\right),
\end{equation}
where $X_{\mK}$ and $Y_{\mK}$ are the unique positive semidefinite 
solutions to the Lyapunov equations below 
\vspace{-1mm}
\begin{subequations}
\begin{align}
X_{\mK} &= \begin{bmatrix} A &  BC_{\mK} \\ B_{\mK} C & A_{\mK} \end{bmatrix}X_{\mK}\begin{bmatrix} A &  BC_{\mK} \\ B_{\mK} C & A_{\mK} \end{bmatrix}^\transpose +  \begin{bmatrix} W & 0 \\ 0 & B_{\mK}VB_{\mK}^\transpose  \end{bmatrix}, \label{eq:LyapunovX_discrete}
\\
Y_{\mK} &= \begin{bmatrix} A &  BC_{\mK} \\ B_{\mK} C & A_{\mK} \end{bmatrix}^\transpose Y_{\mK}\begin{bmatrix} A &  BC_{\mK} \\ B_{\mK} C & A_{\mK} \end{bmatrix} +   \begin{bmatrix} Q & 0 \\ 0 & C_{\mK}^\transpose R C_{\mK} \end{bmatrix}. \label{eq:LyapunovY_discrete}
\end{align}
\end{subequations}
Note that $X_{\mK}$ and $Y_{\mK}$ are closely related to the controllable and observable Gramians of the closed-loop system \eqref{eq:closed_loop_system}. 

From \eqref{eq:LQG_cost_formulation_discrete}, it is not difficult to see that $J_{\texttt{LQG}}(\mK)$ is a rational function in terms of the policy parameter $\mK$, and it is thus a real analytic on $\mathcal{C}_n$. We summarize this as a fact below.%
\begin{fact}
    The LQG cost $J_{\texttt{LQG}}(\mK)$ in \eqref{eq:LQG_policy-optimization}  is real analytic on $\mathcal{C}_n$, and in particular smooth. 
\end{fact}

\begin{figure}[t]
\centering
    \subfigure[]{\includegraphics[width = 0.35\textwidth]{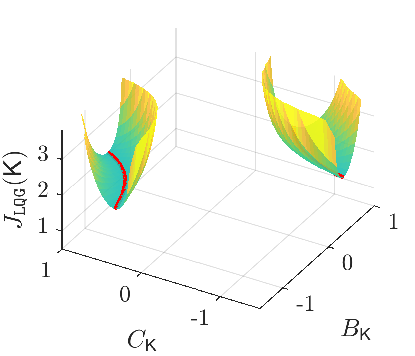}}
    \hspace{15mm}
\subfigure[]{\includegraphics[width=0.35\textwidth]{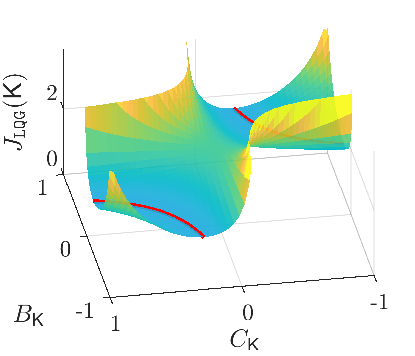}\label{fig:LQG-shape-b}}

    \caption{Non-isolated and disconnected optimal LQG policies (highlighted by the red curves). In both cases, we set $Q = 1$, $R = 1, V = 1, W = 1 $. (a) LQG cost for the open-loop unstable system with $A=1.1, B = C =1$, where we fixed $A_\mK = -0.3069$; (b) another LQG cost for the open-loop stable system with $A=0.9, B = C =1$, where we fixed $A_\mK = -0.1753$. In (b), the origin $B_\mK = 0, C_\mK = 0$ with $A_{\mK} = -0.1753$ is a saddle point.}
    \label{fig:LQG-shape}
\end{figure}

It is also easy to identify examples to confirm the non-convexity of LQG policy optimization \eqref{eq:LQG_policy-optimization} (note that the domain $\mathcal{C}_n$ is already non-convex; see \Cref{fig:feasible_region}). Unlike the LQR, it is known that the problem of LQG problem \eqref{eq:LQG_policy-optimization} has non-unique and non-isolated globally optimal policies in the state-space form $\mathcal{C}_n$. This is not difficult to see since any similarity transformation on one optimal LQG policy $\mK^\star$ leads to another optimal solution that achieves the same cost, i.e., 
\begin{equation} \label{eq:similarity-transformation-LQG}
J_{\texttt{LQG}}(\mK^\star) = J_{\texttt{LQG}}(\mathscr{T}_T(\mK^\star)), \quad \forall T \in \mathrm{GL}_n.
\end{equation}
\Cref{fig:LQG-shape} illustrates the non-convex LQG landscape and the non-isolated/disconnected optimal LQG policies. It is also known that the LQG cost  $J_{\texttt{LQG}}$ is not coercive: there might exist sequences of stabilizing policies $\mK_j \in \mathcal{C}_n$ where $\lim_{j \rightarrow \infty} \mK_j \in \partial \mathcal{C}_n$ such that 
$
    \lim_{j \rightarrow \infty} J_{\texttt{LQG}}(\mK_j) < \infty,  
$
and sequences of stabilizing policies $\mK_j \in \mathcal{C}_n$ where $\lim_{j \rightarrow \infty} \|\mK_j\|_F = \infty$ such that 
$
    \lim_{j \rightarrow \infty} J_{\texttt{LQG}}(\mK_j) < \infty. 
$
The latter fact is easy to see from the effect of similarity transformation \eqref{eq:similarity-transformation-LQG} since $J_{\texttt{LQG}}(\mK)$ is constant for policies that are connected by any $T \in \mathrm{GL}_n$; also see \cite[Example 4.1]{zheng2023benign} for the former fact.
A closed-form expression for the gradient of the LQG cost function $\overline{\nabla} J_{\texttt{LQG}}$ can also be obtained; see \cite{zheng2023benign} for details.

As shown in \Cref{fig:LQG-shape}, the set of stationary points $\{\mK \in \mathcal{C}_n \mid \overline{\nabla} J_{\texttt{LQG}}(\mK) = 0\}$ is not isolated and can be disconnected. Furthermore, there may exist strictly suboptimal spurious stationary points for the LQG control \eqref{eq:LQG_policy-optimization}. This fact can be seen from \Cref{fig:LQG-shape-b}, in which the policy $\mK \in \mathcal{C}_1$ with values $A_\mK =   -0.1753$, $B_\mK = 0$, and $C_\mK = 0$ corresponds to a saddle point.  
Indeed, the following result explicitly characterizes a class of saddle points in LQG control \eqref{eq:LQG_policy-optimization} when the plant dynamics are open-loop stable. These stationary points are spurious and suboptimal whenever the globally  optimal LQG policy corresponds to a nonzero transfer function.

\begin{fact}[Saddle points in LQG] \label{theorem:saddle}
    Suppose \eqref{eq:dynamics_discrete} is open-loop stable. Let $A_\mK = \Lambda\in\mathbb{R}^{n\times n}$ be any stable matrix. Then the zero policy $\mK \in \mathcal{C}_n$ with parameters $A_\mK = \Lambda, B_\mK = \mathbf{0}, C_\mK = \mathbf{0}$ is a stationary point of $J_{\texttt{LQG}}(\mK)$ over $\mathcal{C}_n$, and the corresponding Hessian is either indefinite or zero.
\end{fact}

Due to the existence of spurious saddle points, LQG policy optimization \eqref{eq:LQG_policy-optimization} cannot enjoy the gradient dominance property as the LQR case. The gradient dominance property will also fail for the LQG control even when an observer-based policy parameterization is used \citep{mohammadi2021lack}. Note that the policy $\mK$ in Fact \ref{theorem:saddle} corresponds to a zero transfer function, and the policy just produces a zero input. It has been shown that \textit{all spurious stationary points $\mK$ are non-minimal dynamic policies}, i.e., either $(A_\mK, B_\mK)$ is not controllable, or $(C_\mK, A_\mK)$ is not observable or both. Therefore, we have the following result about globally optimal LQG policies.    

\begin{fact}[Global optimality in LQG]\label{thm:global-optimum}
    All stationary points that correspond to controllable and observable policies are globally optimal in the LQG problem \eqref{eq:LQG_policy-optimization}. These globally optimal policies are related to each other by a similarity transformation. 
\end{fact}

The proof is based on closed-form gradient expressions for controllable and observable stationary points $\mK \in \mathcal{C}_n$ (by letting the gradients %
equal to zero), which are shown to be the same as the optimal solution from Riccati equations. This proof strategy was first used in \citep{hyland1984optimal} to derive first-order necessary conditions for optimal reduced-order controllers. It strongly depends on the assumption of minimality, which fails to deal with non-minimal globally optimal policies for LQG control. Beyond minimal policies, a recent extension of global optimality characterization in LQG control is \cite[Theorem 4.2]{zheng2023benign}, which is based on a notion of \textit{non-degenerate} stabilizing policies. This characterization relies on a more general strategy~of~extended convex lifting that captures a suitable convex re-parameterization of LQG control~\eqref{eq:LQG_policy-optimization}.

\subsubsection{Invariances of LQG: a Quotient Geometry} \label{subsection:symmetry-manifold-structures-LQG}
\begin{figure}[t]
    \centering
    \includegraphics[width = .36\textwidth]{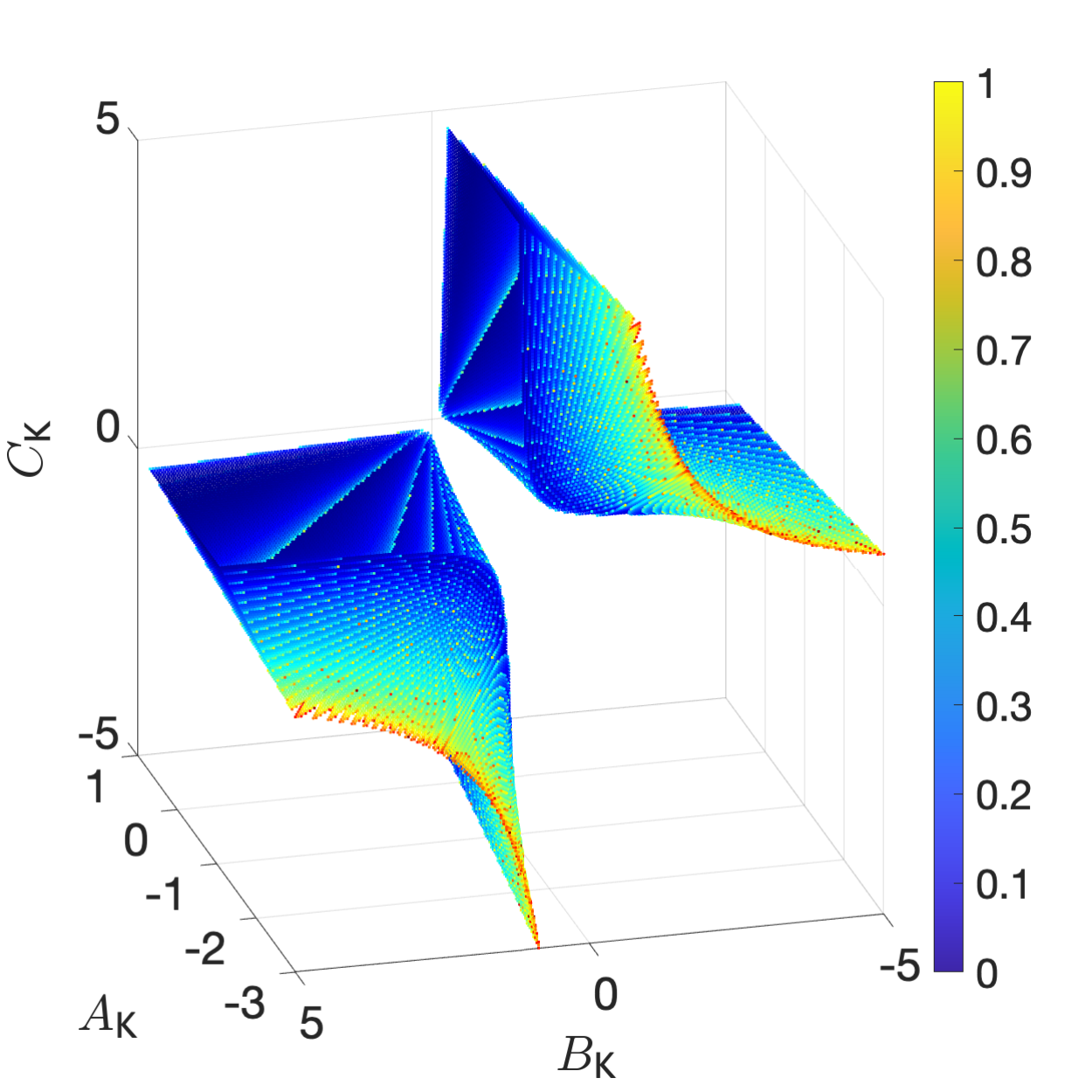}
    \caption{A colored plot of the LQG cost $J_\texttt{LQG}(\cdot)$ over $\mathcal{C}_n$ with $A=1.1, B = C = 1$. Color is determined via an affine transformation of $\log(J_\texttt{LQG}(\mK))$ so that the minimum LQG cost is mapped to 0 and the maximum LQG cost is mapped to 1, over the given plot bounds on $A_\mK, B_{\mK}, C_{\mK}$.}
    \label{fig:LQG_cost_on_manifold_and_orbits}
\end{figure}

This subsection is a sequel to the quotient manifold setup introduced in \S\ref{subsubsec:invariants_of_feedback_policies}. Recall the LQG cost $J_{\texttt{LQG}}$ solely depends on the input-output properties of the closed loop system. If our parameterization of choice is the state-space form in (\ref{eq:internallystabilizing}), then the control design objective will be invariant under similarity transformation: 
$ J_{\texttt{LQG}}(\mK_1) = J_{\texttt{LQG}}(\mK_2),$ for any $\mK_1,\mK_2 \in \mathcal{C}_n$ lying in the same orbit $[\mK_1] = [\mK_2]$. \Cref{fig:LQG_cost_on_manifold_and_orbits} shows an example of such LQG cost on $\mathcal{C}_1$; also compare with the orbits plotted in \Cref{fig:controller_orbit_space}.

The similarity-invariance of LQG causes many hindrances for policy optimization. First, each stationary point lies within an orbit of stationary points, and this orbit will have dimension $n^2$. Therefore, the Hessian of any stationary point will be singular: its nullspace dimension will be at least $n^2$. This may muddle the policy optimization algorithms. Also, local convergence guarantees will be more difficult to obtain since the Hessian at the global minimum cannot be positive-definite. 

Second, since the group product is not coordinate-invariant, we lose the useful invariance property:
\begin{equation*}
    \overline{\nabla} J_{\texttt{LQG}} ( \mathscr{T}_T(\mK))  \not = \mathscr{T}_T(\overline{\nabla} J_{\texttt{LQG}} (\mK)).
\end{equation*} This has impacts on the initialization: if we choose a particularly bad initialization $\|K_0\|_F \gg 0$, then we will always have $\|\overline{\nabla} J_{\texttt{LQG}}(\mK) \|_F \gg 0$. This also implies that gradient descent fails to satisfy the following property known as similarity-equivariance:
\begin{equation}\label{eq:sim-equivariance}
    \mathscr{T}_T(\mK - \alpha \overline{\nabla} J_{\texttt{LQG}}(\mK)) \not = \mathscr{T}_T(\mK) - \alpha \overline{\nabla} J_{\texttt{LQG}}(\mathscr{T}_T(\mK)).
\end{equation} As we will see, we can resolve these issues by introducing a similarity-invariant Riemannian metric and devising a Riemannian gradient descent procedure.

Since $J_\texttt{LQG}(\cdot)$ is similarity-invariant, we can induce a unique cost onto the smooth quotient space $\mathcal{C}_n/\mathrm{GL}_n$ as follows:

\begin{equation*}
    \widetilde{J}_{\texttt{LQG}}([\mK]) := J_{\texttt{LQG}}(\mK')
\end{equation*} where $\mK'$ is any controller $\mK' \in [\mK]$. Now let us equip $\mathcal{C}_n^{\min}$ with the KM metric (\ref{eq:KM-metric}) and $\mathcal{C}_n^{\min}/\mathrm{GL}_n$ with the induced quotient metric. Then the Riemannian gradient is similarity-equivariant:
\begin{equation*}
    \grad J_\texttt{LQG}(\mathscr{T}_T(\mK)) = \mathscr{T}_T(\grad J_\texttt{LQG}(\mK))
\end{equation*} Using the Euclidean retraction, which is similarity-equivariant, we therefore have the fact that a single step of Riemannian gradient descent is similarity equivariant; that is,
\begin{equation*}
    \mathscr{T}_T(\mK - \alpha \grad J_\texttt{LQG}(\mK)) = \mathscr{T}_T(\mK) - \alpha \grad J_\texttt{LQG}(\mathscr{T}_T(\mK)).
\end{equation*} This is a remarkable property. If we perform a similarity transformation on the initial controller $\mK_0$, then the resulting sequence of iterates will also be transformed by that same similarity transformation.

This also implies the following. Let $\mK \in \mathcal{C}_n^{\min}$ and consider $\mK^+ := \mK - \alpha \grad J_\texttt{LQG}(\mK)$. Let $x := [\mK] \in \mathcal{C}_n^{\min}/\mathrm{GL}_n$ and define $x^+ := \mathcal{R}_{x}(-\alpha \grad \widetilde{J}_\texttt{LQG}(x))$. Then $x^+ = [\mK^+]$. In other words, Riemannian gradient descent over the Riemannian quotient manifold coincides with Riemannian gradient descent over the total manifold. We emphasize that these properties hold only when the Riemannian metric is similarity-invariant (\ref{eq:metric-inv}), such as the KM metric, and the retraction is similarity-\textit{equivariant}, such as the Euclidean retraction.

We will end this section with a lemma that shows how this setup reduces the nullspace dimension of Riemannian Hessian at stationary points. The following result quantifies the number of eigenvalues of the Riemannian Hessian with positive, zero, and negative signs---referred to as its ``signature.''
\begin{fact}\label{lem:signature}
    Let $\mathcal{M}$ be a manifold under group action $G$ equipped with a $G$-invariant Riemannian metric. Suppose $\mathcal{M}/G$ is a smooth quotient manifold and consider any $G$-invariant smooth function $J:\mathcal{M} \to \mathbb{R}$. If  $(n_-,n_0,n_+)$ is the signature of $\hess J(x^*)$ at any stationary point $x^* \in \mathcal{M}$ of $J$, then the signature of $\hess \widetilde{J}([x^*])$ is $(n_-,n_0-\dim(G),n_+)$, where $\widetilde{J}([x]):=J(x)$.
\end{fact}
To illustrate a consequence of this lemma, suppose $\mK^*$ is a global minimum of $J_{\texttt{LQG}}$. This implies $\overline{\nabla}^2 J_{\texttt{LQG}}(\mK^*)$ has at least $q^2$ zero eigenvalues. Then, due to the lack of positive-definiteness of the Hessian, iterative updates may not be guaranteed with a linear rate of convergence. But, by this shows that $\hess \widetilde{J}_{\texttt{LQG}}([\mK^*])$ has $q^2$ less zeros. So, if $\dim \ker \overline{\nabla}^2 J_{\texttt{LQG}}(\mK^*) = q^2$ then $\hess \widetilde{J}_{\texttt{LQG}}([\mK^*]) > 0$ which is the property that enables a local linear rate of convergence guaranteed in \cite[Thm. 5.2]{kraisler2024output} for the continuous system dynamics.

\subsection{$\mathcal{H}_\infty$-norm: Systems with Adversarial Noise}

In this section, we return to the state-feedback $\mathcal{H}_\infty$ optimal control problem introduced in \S\ref{subsubsect:state-feedback-robust-control}. 
Setting the initial state $x_0 = 0$ and considering a stabilizing policy $K \in \mathcal{S}$, the $\mathcal{H}_\infty$ performance $J_\infty(\cdot)$ coincides with the square of the $\mathcal{H}_\infty$ norm of the closed-loop transfer function from $w_t$ to a performance measure 
$
z_t := \begin{bmatrix}
    (Q^{1/2} x_t)^\transpose & (R^{1/2} u_t)^\transpose
\end{bmatrix}^\transpose.\
$
Explicitly,
\begin{equation} \label{eq:closed-loop-hinf-norm}
\begin{aligned}
    J_\infty(K) &= \sup_{\|w\|_{l_2} \neq 0} \frac{\|z\|_{l_2}^2}{\|w\|_{l_2}^2}= \left\|\begin{bmatrix}
    Q^{1/2} \\ R^{1/2}K
\end{bmatrix} (zI - A - BK)^{-1}\right\|_\infty^2 \\
&= \sup_{\omega \in [0, 2\pi]}\, \lambda_{\max}\left( (e^{-j\omega} -A - BK)^{-\transpose} (Q + K^\transpose R K) (e^{j\omega} -A - BK)^{-1}\right), 
\end{aligned}
\end{equation}
where $\lambda_{\max}(\cdot)$ denotes the maximal eigenvalue of an Hermitian matrix. %

Unlike the case of LQR or LQG, we usually do not have a closed-form expression to evaluate the $\mathcal{H}_\infty$ norm \eqref{eq:closed-loop-hinf-norm}. One can also use the celebrated KYP lemma to evaluate  it using a linear matrix inequality (LMI), but the solution can take more computational efforts than its counterparts in LQR \eqref{eq:LQR-cost-lyapunov} or LQG \eqref{eq:LQG_cost_formulation_discrete} that only involve solving (linear) Lyapunov equations. %
Classical control techniques typically re-parameterize the non-convex $\mathcal{H}_\infty$ control \eqref{eq:Hinf-policy-optimization} into a convex LMI via a change of variables \citep{boyd1994linear}, or directly characterize a suboptimal solution via solving a single (quadratic) Riccati equation \citep{zhou1988algebraic}. Here, we discuss some geometric aspects of the $\mathcal{H}_\infty$ policy optimization \eqref{eq:Hinf-policy-optimization}.   

\begin{figure}
    \centering
   \subfigure[]{ \includegraphics[width = 0.25\textwidth]{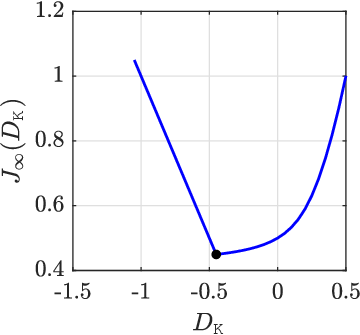}}
    \hspace{15mm}
\subfigure[]{\includegraphics[width=0.33\textwidth]{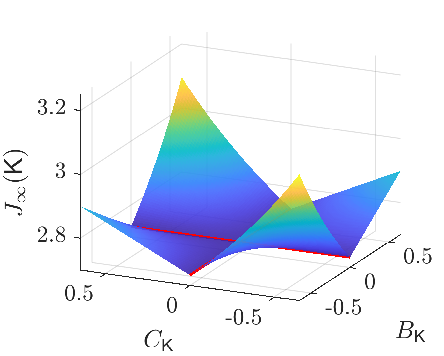}}
    \caption{Non-convexity and non-smoothness of the $\mathcal{H}_\infty$ cost function. (a) illustrates a static output feedback $\mathcal{H}_\infty$ control instance with $u_t = D_{\mK}y_t$ which shows one non-smooth point (highlighted by the black point); (b) corresponds to a dynamic $\mathcal{H}_\infty$ control instance which exhibits a set of nonsmooth points (highlighted by the red lines); see \cite[Example 5.1]{zheng2023benign} for further details. }
    \label{fig:Hinf-landscape} 
\end{figure}

Compared with the LQR and LQG, one major difference in $\mathcal{H}_\infty$ control is that the function $J_{\infty}(K)$ in \eqref{eq:closed-loop-hinf-norm} is non-smooth, i.e., the cost function $J_{\infty}(K)$ may not be differentiable at some feasible points $K \in \mathcal{S}$ (see \Cref{fig:Hinf-landscape} for an illustration). This fact is actually not difficult to see due to two possible sources of non-smoothness in \eqref{eq:closed-loop-hinf-norm}: One from taking the largest eigenvalue of complex matrices, and the other from maximization over  $\omega \in [0,2\pi]$. Indeed, robust control problems were one of the early motivations and applications for non-smooth optimization \citep{lewis2007nonsmooth}. Despite the non-smoothness, the $\mathcal{H}_\infty$ cost function is \textit{locally Lipschitz} and thus differentiable almost everywhere (\Cref{fig:Hinf-landscape} also illustrates this). Thus, we can define the Clarke directional derivative and Clarke subdifferential of $J_{\infty}(K)$ at each feasible policy $K \in \mathcal{S}$.  

Furthermore, the $\mathcal{H}_\infty$ cost function $J_{\infty}(K)$ is known to be ``subdifferentially regular'' in the sense of Clarke \citep{clarke1990optimization} (i.e., the ordinary directional derivative exists and coincides with the Clarke directional derivative for all directions).  
Also, it is known that the discrete-time state-feedback $\mathcal{H}_\infty$ cost function \eqref{eq:closed-loop-hinf-norm} is coercive.
We summarize these properties property below, which to some extend are analogous to those of LQR cost in \Cref{lem:coercive}.

\begin{lemma}\label[lemma]{lemma:H_inf_some_local_Lipschitz}
Suppose $Q \succ 0, R \succ 0$ and $(A,B)$ is stabilizable. Then, the $\mathcal{H}_\infty$ cost function $J_{\infty}\colon \stableK \mapsto \mathbb{R}$, defined in \eqref{eq:closed-loop-hinf-norm}, 
\begin{itemize}
    \item[(a)] is locally Lipschitz, and thus almost everywhere differentiable;
    
    \item[(b)] is subdifferentially regular over the set of stabilizing policies $\mathcal{S}$;
    
    \item[(c)] is coercive\footnote{We remark that this coerciveness property fails to hold in the continuous-time state-feedback $\mathcal{H}_\infty$ control, even when $Q \succ 0, R \succ 0$; see \cite[Fact 4.1]{zheng2024benign}}: $K \to \partial\mathcal{S} \text{ or } \|K\| \to \infty \;\text{ each implies }\; J_{\infty}(K) \to \infty$;
 
    \item[(d)] has compact, path-connected sublevel sets $\mathcal{S}_\gamma = \{K \in \mathcal{S} \mid J_{\infty}(K) \leq \gamma\}$
    for any $\gamma \geq \gamma^\star : = \min_{K} J_{\infty}(K)$.
     
\end{itemize}
\end{lemma}
  
The proof idea for the first two properties is to view $J_{\infty}(K)$ as a composition of a convex mapping $\|\cdot\|_\infty$ and the mapping from  $K$ to a stable closed-loop transfer function that is continuously differentiable over $\mathcal{S}$. 
In the discrete time, the coerciveness of $J_\infty$ can be proved using the positive definiteness of $Q$ and $R$.  
The compactness of sublevel sets follows directly by coercivity. The sublevel set $\mathcal{S}_\gamma$ is in general non-convex but always path-connected.
One can also compute the set of subdifferential  $\partial J_{\infty}(K)$ at each feasible policy $K \in \mathcal{S}$; however, the computation is much more complicated than the smooth LQR or LQG case. We refer the interested reader to \cite[Lemma 5.2]{zheng2023benign} and \cite{apkarian2006nonsmooth} for more details. %
Despite the non-convex and non-smoothness, we have a global optimality characterization for \eqref{eq:Hinf-policy-optimization}. 

\begin{fact}[Global optimality in $\mathcal{H}_\infty$ control] \label{theorem:hinf-global-optimality}
    Consider the state-feedback $\mathcal{H}_\infty$ control \eqref{eq:Hinf-policy-optimization} with $Q \succ 0, R \succ 0$. Any Clarke stationary point is globally optimal.  
\end{fact}

The high-level proof of the above results proceed as follows:
it is known that \eqref{eq:Hinf-policy-optimization} admits an equivalent convex reformulation by a change of variables, and this change of variables can be designed as a diffeomorphism between non-convex policy optimization and its convex reformulation; this diffeomorphism then allows us to certify global optimality in original non-smooth and non-convex $\mathcal{H}_\infty$ control; see \cite[Theorem 1]{guo2022global} and \cite[Corollary 4.1]{zheng2024benign} for details. This idea has further been characterized into a framework of extended convex lifting (ECL) in \citep{zheng2023benign,zheng2024benign}, which bridges the gap between
non-convex policy optimization and convex reformulations in a range of control problems. %

\section{Algorithmic Implications}  \label{data-driven}

In the context of reinforcement learning and control, geometric perspectives on policy optimization facilitate the development of data-driven algorithms that can emulate various first- and second-order policy iteration schemes. These approaches typically involve synthesizing a first- or second-order oracle using available performance measure information. This framework provides a basis for comparing the data efficiency of different techniques in terms of sample complexity. Specifically, it addresses how many function calls to the oracle are necessary to achieve the desired level of optimality %
when the algorithm has access only to the oracle rather than the explicit problem parameters.

\subsection{Convergence of Policy Optimization Algorithms} \label{section:exact-gradient}

Despite the non-convexity of the LQR problem in the policy parameters $K$, the analysis of the domain manifold and the properties of LQR cost in \S\ref{susec:standard-LQR-perf}, in particular the gradient dominance property has enabled establishing the following global linear convergence guarantees of gradient descent algorithms \citep{fazel2018global,bu_lqr_2019}. 
This linear convergence result is mainly due to the \textit{coerciveness}, \textit{smoothness} over any sublevel set, and \textit{gradient dominance} of the LQR cost function $J_{\texttt{LQR}}(K)$ (see \Cref{lem:coercive}).

\begin{fact}
    Starting from any feasible $K_0 \in \stableK$, a  small enough (but constant) step-size $\eta$ remains stabilizing for the gradient descent updates $K^+ = K -\eta \overline{\nabla} J_{\texttt{LQR}}(K)$ which converges to the optimal LQR policy $K^*$ at a linear rate.
\end{fact}

As discussed in \Cref{rem:hewer}, the algebraic update of Hewer's algorithm can be described as a ``Riemannian quasi-Newton'' update and we can provide an alternative proof for its global convergence \cite{talebi_policy_2023}.

\begin{fact}
    Starting from any feasible $K_0 \in \stableK$, the unit step size remains stabilizing for the Hewer's update $K^+ = K + \widehat{V}$ with $\widehat{V}$ solving $\widehat{H}_{K}\widehat{V} = - \grad{J}_{\texttt{LQR}}(K)$ which converges to the optimal LQR policy $K^*$ at a quadratic rate.
\end{fact}

In the specific context of Hewer's update, the input-output system trajectory can be directly utilized to obtain a positive definite approximation of the Riemannian Hessian \( \hat{H}_K \) and the Riemannian gradient \( \grad J_{\texttt{LQR}(K)} \) through a recursive least squares scheme \cite{bradtke_adaptive_1994}. This approach can be extended to solve constrained LQR problems, as demonstrated in \cite{alemzadeh_data-driven_2024}, which focuses on learning policies that adhere to a communication/information graph in a large network of homogeneous systems. Extensions to any linearly constrained policies, including static output feedback and structured LQR problems, utilizing the stability certificate idea in  \Cref{lem:stability-cert} are explored in \cite{talebi_policy_2023} and summarized below:

\begin{fact}
    Starting from any feasible $K_0 \in \substableK = \stableK\cap \constraint$, the stepsize $\eta_K = \min(1,s_K)$ remains stabilizing for the Riemannian Newton update  $K^+ = K + \eta_K V_K$ with $V_K$ solving $\hess J_{\texttt{LQR}}|_K(V_K) = -\grad J_{\texttt{LQR}}(K)$ and its variants (by replacing $\hess$ with $\euchess$). Furthermore, any non-degenerate local minimum is contained in a neighborhood on which the generated sequence of polices remains therein and converges to the local minimum fast--at a linear rate that eventually becomes quadratic.
\end{fact}
\noindent Inspired by this result, recently an online optimistic version of these updates are studied in \cite{chang_regret_2024} with regret bound guarantees.

Let us move to the policy optimization for LQG control \eqref{eq:LQG_policy-optimization} over the dynamic output-feedback policies $\mathcal{C}_n$. Despite being non-convex, the geometrical properties of $\mathcal{C}_n$ 
and the LQG cost function $J_{\texttt{LQG}}(\cdot)$ (in \S\ref{subsec:LQG}) can ensure some favorable properties of policy gradient algorithms. While the feasible region $\mathcal{C}_n$ can be path-connected, Fact \ref{Fact:disconnectivity} ensures that the two path-connected components are identical from an input-output perspective. Thus, when applying policy search algorithms to solve LQG problem \eqref{eq:LQG_policy-optimization}, it makes no difference to search over either path-connected component. In addition, if a sequence of gradient iterates converges to a point, Fact \ref{thm:global-optimum} further allows us to verify whether the limit point is a globally optimal solution to the LQG control \eqref{eq:LQG_policy-optimization}. The following is an immediate corollary of this fact. 

\begin{fact} \label{fact:Gradient_Descent-LQG}
    Consider a gradient descent algorithm $\mK_{t+1} = \mK_{t} - \alpha_t \overline{\nabla} J_{\texttt{LQG}}(\mK_t)$ for the LQG problem~\eqref{eq:LQG_policy-optimization}, where $\alpha_t$ is a step size.  Suppose $\inf_t\alpha_t>0$ and the iterates $\mK_{t}$ converge to a point $\mK^*$. Then $\mK^*$ is globally optimal if it is a minimal controller. 
\end{fact}

We emphasize that there are two major limitations in Fact \ref{fact:Gradient_Descent-LQG}: 1) it does not address the case when the limit point has a vanishing gradient ($\overline{\nabla} J_{\texttt{LQG}}(\mK^*) = 0$) but is non-minimal; 2) it does not offer conditions to ensure convergence of the gradient descent iterates.  Indeed, Fact \ref{theorem:saddle} has revealed that there may exist strictly suboptimal saddle points for non-minimal LQG policies in~\eqref{eq:LQG_policy-optimization}. If a stationary point does not correspond to a minimal (aka, controllable and observable) policy, we cannot confirm optimality.  Furthermore, some saddle points are \textit{high-order} in the sense that they have degenerate Hessian (e.g., the corresponding Hessian is zero), and thus there is no escaping direction that is a key element in the developments on perturbed gradient methods to avoid saddles \citep{jin2017escape}. Recently, \cite{zheng2022escaping} introduced a new perturbed policy gradient (PGD) algorithm to escape a class of spurious stationary points (including high-order saddles). One key idea is to use a novel reparameterization procedure that converts the iterate from a high-order saddle to a strict saddle, from which standard random perturbations in gradient descent can escape efficiently. 

The inherent symmetry induced by similarity transformation still makes the convergence conditions of ordinary gradient descent methods hard to derive. With the Riemannian quotient manifold setup $\mathcal{C}_n^{\min}/\mathrm{GL}_n$ in \S\ref{subsection:symmetry-manifold-structures-LQG}, a local linear convergence rate for Riemannian gradient methods is derived in \cite{kraisler2024output}--whenever the iterates are close to a globally optimal policy. We summarize this result as follows.

\begin{fact}
    Let $\mathcal{C}_n^{\min}/\mathrm{GL}_n$ be the orbit space of full-order minimal dynamic feedback controllers (for the continuous LTI dynamics) modulo similarity transformation. %
    Under certain regularity conditions on $J_{\texttt{LQG}}$ \cite[Assumption 5.1]{kraisler2024output}, let $\mK^*$ be an optimal LQG controller. Consider the Riemannian gradient descent updates under the KM metric $\langle .,. \rangle^{\mathrm{KM}}$ and the Euclidean retraction, written as $\mK_{t+1} = \mK_t - \alpha \grad J_{\texttt{LQG}}(\mK_t)$ with a sufficiently small step size $\alpha > 0$. Then, there exists a neighborhood of $[\mK^*]$ in which if we initialize $K_0 \in [\mK^*]$, then  $\lim_{t \to \infty} [\mK_t] = [\mK^*]$. That is, the orbit of $[\mK]$ converges to the orbit of $[\mK^*]$. %
    Furthermore, the rate of convergence is linear.
\end{fact}

\subsection{Oracle-based Data-driven Algorithms}

The discussions in \S\ref{section:exact-gradient} require exact (Riemannian) gradient information. In model-free scenarios, it is possible to evaluate the performance measure $ J(\theta) $ for a given set of policy parameters $ \theta $ that determines an input $ \mathbf{u} = \pi_\theta(\mathbf{x}) $. As will become clear in the following discussion, it is reasonable to construct oracles based on these function evaluations, which may differ depending on the specific method used for gradient estimation. This approximate evaluation of $ J(\theta) $ is feasible whenever its explicit form is known and the system's input-output trajectory $(\mathbf{u}, \mathbf{y})$ can be obtained through independent experiments. Like other sample-based techniques, the performance of each approximation can be quantified by evaluating the bias-variance trade-off, which directly impacts the probabilistic convergence guarantees of these methods.

Probably the most natural of these approaches is the Finite Difference Method, where the gradient is estimated by relating it back to the performance difference of randomly selected perturbations in the \( d \)-dimensional policy parameters \( \theta \in \Theta \) through smoothing techniques. In particular, one may approximate \( J(\theta) \) by the following averaging/smoothing:
\[J(\theta) \approx \hat{J}(\theta) \coloneqq \mathbb{E}_{\nu \sim \mathrm{Uni}\{\mathbb{B}^d\}} J(\theta + \varepsilon \nu),\]
where \( \mathbb{B}^d \) denotes the unit \( d \)-dimensional ball and \( \varepsilon \) is a small radius. Then, the gradient can be approximated by \cite[Lemma 1]{flaxman_online_2004}:
\[\overline{\nabla} \hat{J}(\theta) = \mathbb{E}_{U \sim \mathrm{Uni}\{\mathbb{S}^d\}} \left[\frac{J(\theta + \varepsilon U)}{\varepsilon/d} U \right],\]
where \( \mathbb{S}^d \) denotes the unit \( d \)-dimensional sphere.
In line with this, the so-called two-point approximation from finite samples can be expressed as \cite{spall_overview_1998}:
\[
\overline{\nabla} \hat{J}(\theta) = \frac{1}{N}\sum_{i=1}^N \frac{J(\theta + \varepsilon U_i) - J(\theta - \varepsilon U_i)}{2\varepsilon/d} U_i,
\]
where \( U_i \) is a (feasible) randomly selected unit vector, \( \varepsilon \) is a small enough perturbation parameter, and \( N \) is the number of samples. Note that if the perturbation size is particularly small, we can ensure the feasibility of the perturbed policy \( \theta + \varepsilon U \), especially when \( \Theta \) has a relatively open structure, as in the case of the set of static/constrained/dynamic stabilizing policies \( \stableK \), \( \substableK \), or \( \mathcal{C}_q \). For data efficiency, the two-point approximation can be further reduced to the so-called one-point approximation as follows \cite{flaxman_online_2004}:
\[\overline{\nabla}_\theta \hat{J}(\theta) \approx \frac{1}{N}\sum_{i=1}^N \frac{J(\theta + \varepsilon U_i)}{\varepsilon/d} U_i.\]
These estimations, for example, enables learning optimal LQR policy form input output trajectories with (probabilistic) global convergence guarantees \cite{fazel2018global}.

Similarly, analogous arguments can be followed to estimate the Hessian of the cost from additional independent samples \cite{shen_hessian_nodate}. However, these techniques often lead to high variance issues, which can be mitigated by introducing an initial state-dependent ``baseline'' \( J(\theta;x_0) \) for approximating these variations \cite{grondman_survey_nodate, kakade_natural_2001}:
\[\overline{\nabla}_\theta \hat{J}(\theta) \approx \frac{1}{N}\sum_{i=1}^N \frac{J(\theta + \varepsilon U_i; x_0^i) - J(\theta; x_0^i)}{\varepsilon/d} U_i,\]
where the baseline is independently approximated by \( J(\theta; x_0^i) = \mathbb{E}_{U \sim \mathrm{Uni}\{\mathbb{S}^d\}} \left[\frac{J(\theta + \varepsilon U)}{\varepsilon/d} U \right] \) for each initial state \( x_0^i \). This technique is adopted, for example, in \cite{takakura_structured_2024} to reducing variance of learning output feedback LQR policy. A similar approach can be applied for Hessian approximation.

It is also worth noting that these data generation procedures for model-free function evaluation, often require \textit{a priori} access to a stabilizing policy for the underlying system dynamics. This relates to online stabilization problems and its intricate geometry from its fundamental limitations \cite{talebi_regularizability_2022} to algorithm design; see for example \cite{yu_online_2023}.

\subsection{Optimal Estimation Problems}
Another recent development has seen the translation of policy optimization techniques, originally developed for optimal control, to optimal estimation problems through the profound ``duality relation'' between these two setups \cite{talebi_data-driven_2023}. In the optimal estimation context, the mean-squared estimation (MSE) error can be naturally expressed as an average
\[J_{\texttt{MSE}}(\theta) = \mathbb{E}_{\mathbf{y}}\mathsf{SE}(\theta,\mathbf{y}),\]
where \( \mathsf{SE}(\theta,\mathbf{y}) \) denotes the squared estimation error. The gradient of this error can be computed for any observed trajectory \( \mathbf{y} \) and given (now called) estimation policy \( \theta \). This gradient, \( \overline{\nabla}_\theta \mathsf{SE}(\theta,\mathbf{y}) \), can be approximated using finite-length output trajectories \( \mathbf{y}_T \). This results in a natural gradient approximation scheme based on \( N \) finite-length output trajectories as follows:
\[
\overline{\nabla} J_{\texttt{MSE}}(\theta) \approx \overline{\nabla} \hat J_T(\theta) = \frac{1}{N}\sum_{i=1}^N \overline{\nabla}_\theta \mathsf{SE}(\theta,\mathbf{y}_T^i).
\]
Finally, an analysis of the bias-variance trade-off enables the establishment of probabilistic guarantees for the convergence of Stochastic Gradient Descent (SGD) for \( \theta \) to the optimal Kalman estimation policy \cite{talebi_data-driven_2023}:

\begin{fact}
    Suppose the system is observable and both dynamic and measurement noise are bounded. Consider the stochastic gradient descent on the estimation policy $\theta^+ = \theta - \eta \overline{\nabla} \hat J_T(\theta)$ with small enough stepsize $\eta$. Then, with high probability, it converges linearly and globally (from any initial stabilizing policy) to the optimal Kalman estimation policy.
\end{fact}

\subsection{Broader Implications: Iterative Learning Procedures}
Other policy parameterization techniques have seen significant success over the last couple of decades, particularly through Linear Matrix Inequality (LMI) techniques, which enable the formulation of stability, robustness, and other performance considerations. These approaches often rely on parameterizing policies in specific ways, such as Youla parameterization, which can be heavily dependent on the underlying system model.

However, these ``model-dependent formulations'', such as those involving Riccati equations and LMI techniques, have limitations when it comes to generalizing across nonlinear dynamics and complex policy parameterizations. In contrast, the complete policy optimization approach offers greater generalization power, particularly for nonlinear dynamics and policies, such as those using neural networks. Additionally, it simplifies the imposition of direct constraints on the synthesized input signal. Incorporating such constraints within those model-dependent frameworks is not straightforward, making the complete policy optimization approach more versatile and robust for a broader range of applications.

\section{Summary and Outlook} \label{summary}

In this survey, we have provided an overview of recent progress on 
understanding geometry of policy optimization and its algorithmic implications. This has been pursued both in terms of the
static and dynamic stabilization problems, as well as the how
control performance objectives interact with this set.
The implications of such geometric perspective on policy optimization
for developing first order methods for control design, 
as well as their model-free data driven realizations 
are also discussed.

Some of key ideas that underlie our presentation include
developing a geometric machinery to reason about
fundamental (complexity) bounds for feedback design, both in terms of
model parameters as well as available data. 
For example, we advocate that understanding the geometry of the cost, when constrained to submanifolds of stabilizing feedback gains, is crucial for devising efficient model-based and model-free algorithms for robust and optimal designs, be it in terms of homotopies, escaping saddle points, conditioning, or effective use of symmetries.

\section{Notes and Commentary} \label{notes-commentary}

Throughout this manuscript, we reference \citep{lee_introduction_2018} for standard geometric notions such as Riemannian metric, connection, vector field, gradient, Hessian, and Weingarten mapping.
In \S\ref{section:stablizing-policies}, the topological properties of static stabilizing policies are from \citep{bu_topological_2021}, with earlier results in \citep{ohara1992differential,fam1978canonical,ober1987topology}. For static state-feedback Hurwitz stabilizing policies in continuous-time LTI systems, see \citep{ohara_differential_1992}. 
The geometric PO ideas on \ac{slqr} and \ac{olqr} are reviewed from \citep{talebi_riemannian_2022}, and further details on the gradient, Hessian, and Christoffel symbols are in \citep{talebi_policy_2023}.

The results in on dynamic feedback synthesis for LQG in \Cref{subsec:LQG} are based on \citep{tang2023analysis,kraisler2024output}, and other related results can be found in  \citep{Duan2024LQcontrol,hu2022connectivity}. The topological properties of stabilizing dynamic feedback policies are discussed in \citep{tang2023analysis}, with Fact \ref{Fact:disconnectivity} adapted from \cite[Theorems 1 \& 2]{tang2023analysis}. Detailed computations for the examples of stabilizing policies in \Cref{fig:feasible_region} can be found in \cite[Example 11]{tang2023analysis}. 
Quotient spaces of linear systems are first studied by Kalman and Hazelwinkel, known as geometric linear system theory, in the early 1970s \citep{beckmann_invariants_1976,byrnes_fine_1980}. This perspective focuses on the state-space forms of systems and their algebraic-geometric properties under the feedback. Herein, we focus on the space of stabilizing policies and their symmetries under similarity transformations.
The (KM) Riemannian metric is introduced in \citep{kraisler2024output}, differing slightly from the KM metric in \citep{krishnaprasad_families_1983a} and studied in \citep{afsari_bundle_2017}. \Cref{thm:Riem-metric-dynamic} is proved in \cite[Thm. 3.2]{kraisler2024output}. For abstract quotient spaces and smooth quotient manifold conditions, see \cite[\S21]{lee_introduction_2010}.

In \S\ref{section:LQ-performance}, \Cref{lem:stability-cert} is reviewed from \citep[Lemma 4.1]{talebi_policy_2023}. Hewer's algorithm in \Cref{rem:hewer} is introduced in \citep{hewer_iterative_1971}. 
The properties of $J_\texttt{LQR}$ in \Cref{lem:coercive} are studied in \citep{fazel2018global,bu_lqr_2019,talebi_policy_2023}, with similar properties for Mean-squared error in state estimation in \citep{talebi_data-driven_2023}. 
The results of \S\ref{subsec:LQG} are adapted from \citep{tang2023analysis,kraisler2024output,zheng2022escaping}, with related developments in \citep{hu2022connectivity,zheng2023benign,umenberger2022globally,zhang2023learning}. 
Facts \ref{theorem:saddle} and \ref{thm:global-optimum} are adapted from \cite[Theorem 5 and 6]{tang2023analysis}.
Fact \ref{lem:signature} is formally proved in \cite[Lem 5.3]{kraisler2024output}. For the KYP lemma, see \citep{rantzer1996kalman}, and for $G$-invariant metrics and $G$-equivariance retraction, see \cite[Sect. 9.9]{boumal_introduction_2023}.
 
\Cref{lemma:H_inf_some_local_Lipschitz} aggregates results from various resources. The non-smoothness of $\mathcal{H}_\infty$ cost is in \citep{apkarian2006nonsmooth,apkarian2006nonsmooth2}, with subdifferential regularity from \citep{clarke1990optimization}. Recent discussions are in \cite[Proposition 2]{guo2022global}, \cite[Proposition 1]{tang2023global}, with coercivity in \cite[Lemma 1]{guo2022global} and connectivity of sublevel sets in \cite[Lemma 2]{guo2022global} and \citep{hu2022connectivity}. The subdifferentials of the $\mathcal{H}_\infty$ cost function are computed in \cite[Lemma 5.2]{zheng2023benign}. Fact \ref{theorem:hinf-global-optimality} is first proved in \cite[Theorem 1]{guo2022global} for discrete-time dynamics, with the continuous-time state-feedback $\mathcal{H}_\infty$ control counterpart in \cite[Corollary 4.1]{zheng2024benign}.
At the time of writing, geometrical properties for output-feedback $\mathcal{H}_\infty$ control are under active investigation; see e.g., \citep{zheng2023benign,tang2023global,guo2024complexity}.

In addition to the state-feedback $\mathcal{H}_\infty$ control, some recent studies have also investigated policy optimization in other control problems with robustness features. For example, policy optimization for linear risk-sensitivity control and a general mixed $\mathcal{H}_2$/$\mathcal{H}_\infty$ is studied in \cite{zhang2020policy}, where a notion of implicit regularization is introduced to deal with the challenge of lacking coerciveness in the mixed design. Model-free $\mu$-synthesis was studied in \cite{keivan2022model} and global convergences for risk-constrained LQR are recently investigated in \cite{zhao2023global}. 
In the realm of data-driven policy optimization algorithms, notable methods include REINFORCE \cite{williams_simple_1992} and Proximal Policy Optimization (PPO) \cite{schulman_proximal_2017} which operate within the context of Markov Decision Processes (MDPs), leveraging the finiteness of state and action domains. They estimate the gradient by using the likelihood ratio method \cite{sutton_policy_1999}. 

\section{Acknowledgements}
S. Talebi and N. Li are partially supported by NSF AI institute 2112085. Y. Zheng is supported in part by NSF ECCS-2154650 and NSF CAREER 2340713. S. Kraisler and M. Mesbahi have been supported by NSF grant ECCS-2149470 and AFOSR grant FA9550-20-1-0053.

\bibliographystyle{elsarticle-harv} 
\bibliography{references,references-2}

\begin{thebibliography}{89}
\expandafter\ifx\csname natexlab\endcsname\relax\def\natexlab#1{#1}\fi
\providecommand{\url}[1]{\texttt{#1}}
\providecommand{\href}[2]{#2}
\providecommand{\path}[1]{#1}
\providecommand{\DOIprefix}{doi:}
\providecommand{\ArXivprefix}{arXiv:}
\providecommand{\URLprefix}{URL: }
\providecommand{\Pubmedprefix}{pmid:}
\providecommand{\doi}[1]{\href{http://dx.doi.org/#1}{\path{#1}}}
\providecommand{\Pubmed}[1]{\href{pmid:#1}{\path{#1}}}
\providecommand{\bibinfo}[2]{#2}
\ifx\xfnm\relax \def\xfnm[#1]{\unskip,\space#1}\fi
\bibitem[{Ackermann(1980)}]{ackermann_parameter_1980}
\bibinfo{author}{Ackermann, J.}, \bibinfo{year}{1980}.
\newblock \bibinfo{title}{Parameter space design of robust control systems}.
\newblock \bibinfo{journal}{IEEE Transactions on Automatic Control}
  \bibinfo{volume}{25}, \bibinfo{pages}{1058--1072}.
\newblock \DOIprefix\doi{10.1109/TAC.1980.1102505}.
\bibitem[{Afsari and Vidal(2017)}]{afsari_bundle_2017}
\bibinfo{author}{Afsari, B.}, \bibinfo{author}{Vidal, R.},
  \bibinfo{year}{2017}.
\newblock \bibinfo{title}{Bundle {Reduction} and the {Alignment} {Distance} on
  {Spaces} of {State}-{Space} {LTI} {Systems}}.
\newblock \bibinfo{journal}{IEEE Transactions on Automatic Control}
  \bibinfo{volume}{62}, \bibinfo{pages}{3804--3819}.
\newblock \URLprefix \url{http://ieeexplore.ieee.org/document/7872404/},
  \DOIprefix\doi{10.1109/TAC.2017.2678839}.
\bibitem[{Agarwal et~al.(2021)Agarwal, Kakade, Lee and
  Mahajan}]{Agarwal2021-ss}
\bibinfo{author}{Agarwal, A.}, \bibinfo{author}{Kakade, S.M.},
  \bibinfo{author}{Lee, J.D.}, \bibinfo{author}{Mahajan, G.},
  \bibinfo{year}{2021}.
\newblock \bibinfo{title}{On the theory of policy gradient methods: optimality,
  approximation, and distribution shift}.
\newblock \bibinfo{journal}{J. Mach. Learn. Res.} \bibinfo{volume}{22},
  \bibinfo{pages}{4431--4506}.
\bibitem[{Alemzadeh et~al.(2024)Alemzadeh, Talebi and
  Mesbahi}]{alemzadeh_data-driven_2024}
\bibinfo{author}{Alemzadeh, S.}, \bibinfo{author}{Talebi, S.},
  \bibinfo{author}{Mesbahi, M.}, \bibinfo{year}{2024}.
\newblock \bibinfo{title}{Data-{Driven} {Structured} {Policy} {Iteration} for
  {Homogeneous} {Distributed} {Systems}}.
\newblock \bibinfo{journal}{IEEE Transactions on Automatic Control} ,
  \bibinfo{pages}{1--15}\URLprefix
  \url{https://ieeexplore.ieee.org/abstract/document/10436325},
  \DOIprefix\doi{10.1109/TAC.2024.3366038}. \bibinfo{note}{conference Name:
  IEEE Transactions on Automatic Control}.
\bibitem[{Amari(1998)}]{amari_natural_1998}
\bibinfo{author}{Amari, S.i.}, \bibinfo{year}{1998}.
\newblock \bibinfo{title}{Natural {Gradient} {Works} {Efficiently} in
  {Learning}}.
\newblock \bibinfo{journal}{Neural Computation} \bibinfo{volume}{10},
  \bibinfo{pages}{251--276}.
\newblock \URLprefix \url{https://doi.org/10.1162/089976698300017746},
  \DOIprefix\doi{10.1162/089976698300017746}.
\bibitem[{Amari(2016)}]{Amari2016-do}
\bibinfo{author}{Amari, S.I.}, \bibinfo{year}{2016}.
\newblock \bibinfo{title}{Information Geometry and Its Applications}.
\newblock \bibinfo{publisher}{Springer}.
\bibitem[{Anderson and Moore(2007)}]{anderson_optimal_2007}
\bibinfo{author}{Anderson, B.D.}, \bibinfo{author}{Moore, J.B.},
  \bibinfo{year}{2007}.
\newblock \bibinfo{title}{Optimal control: linear quadratic methods}.
\newblock \bibinfo{publisher}{Courier Corporation}.
\bibitem[{Apkarian and Noll(2006a)}]{apkarian2006nonsmooth}
\bibinfo{author}{Apkarian, P.}, \bibinfo{author}{Noll, D.},
  \bibinfo{year}{2006}a.
\newblock \bibinfo{title}{Nonsmooth {$H_\infty$} synthesis}.
\newblock \bibinfo{journal}{IEEE Transactions on Automatic Control}
  \bibinfo{volume}{51}, \bibinfo{pages}{71--86}.
\bibitem[{Apkarian and Noll(2006b)}]{apkarian2006nonsmooth2}
\bibinfo{author}{Apkarian, P.}, \bibinfo{author}{Noll, D.},
  \bibinfo{year}{2006}b.
\newblock \bibinfo{title}{Nonsmooth optimization for multidisk {$H_\infty$}
  synthesis}.
\newblock \bibinfo{journal}{European Journal of Control} \bibinfo{volume}{12},
  \bibinfo{pages}{229--244}.
\bibitem[{{\AA}str{\"o}m and Murray(2021)}]{aastrom2021feedback}
\bibinfo{author}{{\AA}str{\"o}m, K.J.}, \bibinfo{author}{Murray, R.},
  \bibinfo{year}{2021}.
\newblock \bibinfo{title}{Feedback systems: an introduction for scientists and
  engineers}.
\newblock \bibinfo{publisher}{Princeton university press}.
\bibitem[{Bertsekas(2012)}]{Bertsekas2012-oq}
\bibinfo{author}{Bertsekas, D.}, \bibinfo{year}{2012}.
\newblock \bibinfo{title}{Dynamic Programming and Optimal Control: Volume {I}}.
\newblock \bibinfo{publisher}{Athena Scientific}.
\bibitem[{Bertsekas(2011)}]{Bertsekas2011-hi}
\bibinfo{author}{Bertsekas, D.P.}, \bibinfo{year}{2011}.
\newblock \bibinfo{title}{Approximate policy iteration: a survey and some new
  methods}.
\newblock \bibinfo{journal}{J. Control Theory Appl.} \bibinfo{volume}{9},
  \bibinfo{pages}{310--335}.
\bibitem[{Bertsekas(2017)}]{Bertsekas2017-fi}
\bibinfo{author}{Bertsekas, D.P.}, \bibinfo{year}{2017}.
\newblock \bibinfo{title}{Value and policy iterations in optimal control and
  adaptive dynamic programming}.
\newblock \bibinfo{journal}{IEEE Trans. Neural Netw. Learn. Syst.}
  \bibinfo{volume}{28}, \bibinfo{pages}{500--509}.
\bibitem[{Boumal(2023)}]{boumal_introduction_2023}
\bibinfo{author}{Boumal, N.}, \bibinfo{year}{2023}.
\newblock \bibinfo{title}{An introduction to optimization on smooth manifolds}.
\newblock \bibinfo{publisher}{Cambridge University Press},
  \bibinfo{address}{Cambridge ; New York, NY}.
\bibitem[{Boyd et~al.(1994)Boyd, El~Ghaoui, Feron and
  Balakrishnan}]{boyd1994linear}
\bibinfo{author}{Boyd, S.}, \bibinfo{author}{El~Ghaoui, L.},
  \bibinfo{author}{Feron, E.}, \bibinfo{author}{Balakrishnan, V.},
  \bibinfo{year}{1994}.
\newblock \bibinfo{title}{Linear Matrix Inequalities in System and Control
  Theory}.
\newblock \bibinfo{publisher}{Society for Industrial and Applied Mathematics}.
\bibitem[{Bradtke et~al.(1994)Bradtke, Ydstie and
  Barto}]{bradtke_adaptive_1994}
\bibinfo{author}{Bradtke, S.}, \bibinfo{author}{Ydstie, B.},
  \bibinfo{author}{Barto, A.}, \bibinfo{year}{1994}.
\newblock \bibinfo{title}{Adaptive linear quadratic control using policy
  iteration}, in: \bibinfo{booktitle}{Proceedings of 1994 {American} {Control}
  {Conference} - {ACC} '94}, \bibinfo{publisher}{IEEE},
  \bibinfo{address}{Baltimore, MD, USA}. pp. \bibinfo{pages}{3475--3479}.
\newblock \URLprefix \url{http://ieeexplore.ieee.org/document/735224/},
  \DOIprefix\doi{10.1109/ACC.1994.735224}.
\bibitem[{Brockett(2014)}]{Brockett2014-zg}
\bibinfo{author}{Brockett, R.}, \bibinfo{year}{2014}.
\newblock \bibinfo{title}{The early days of geometric nonlinear control}.
\newblock \bibinfo{journal}{Automatica} \bibinfo{volume}{50},
  \bibinfo{pages}{2203--2224}.
\bibitem[{Bu et~al.(2019)Bu, Mesbahi, Fazel and Mesbahi}]{bu_lqr_2019}
\bibinfo{author}{Bu, J.}, \bibinfo{author}{Mesbahi, A.},
  \bibinfo{author}{Fazel, M.}, \bibinfo{author}{Mesbahi, M.},
  \bibinfo{year}{2019}.
\newblock \bibinfo{title}{{LQR} through the {Lens} of {First} {Order}
  {Methods}: {Discrete}-time {Case}}.
\newblock \URLprefix \url{http://arxiv.org/abs/1907.08921},
  \DOIprefix\doi{10.48550/arXiv.1907.08921}. \bibinfo{note}{arXiv:1907.08921
  [cs, eess, math]}.
\bibitem[{Bu et~al.(2021)Bu, Mesbahi and Mesbahi}]{bu_topological_2021}
\bibinfo{author}{Bu, J.}, \bibinfo{author}{Mesbahi, A.},
  \bibinfo{author}{Mesbahi, M.}, \bibinfo{year}{2021}.
\newblock \bibinfo{title}{On {Topological} {Properties} of the {Set} of
  {Stabilizing} {Feedback} {Gains}}.
\newblock \bibinfo{journal}{IEEE Transactions on Automatic Control}
  \bibinfo{volume}{66}, \bibinfo{pages}{730--744}.
\newblock \URLprefix
  \url{https://ieeexplore.ieee.org/abstract/document/9086139},
  \DOIprefix\doi{10.1109/TAC.2020.2992510}. \bibinfo{note}{conference Name:
  IEEE Transactions on Automatic Control}.
\bibitem[{Chang and Shahrampour(2024)}]{chang_regret_2024}
\bibinfo{author}{Chang, T.J.}, \bibinfo{author}{Shahrampour, S.},
  \bibinfo{year}{2024}.
\newblock \bibinfo{title}{Regret {Analysis} of {Policy} {Optimization} over
  {Submanifolds} for {Linearly} {Constrained} {Online} {LQG}}.
\newblock \URLprefix \url{http://arxiv.org/abs/2403.08553}.
  \bibinfo{note}{arXiv:2403.08553 [cs, eess, math]}.
\bibitem[{Clarke(1990)}]{clarke1990optimization}
\bibinfo{author}{Clarke, F.H.}, \bibinfo{year}{1990}.
\newblock \bibinfo{title}{Optimization and Nonsmooth Analysis}.
\newblock \bibinfo{publisher}{Society for Industrial and Applied Mathematics}.
\bibitem[{Doyle et~al.(2013)Doyle, Francis and Tannenbaum}]{Doyle2013-vh}
\bibinfo{author}{Doyle, J.C.}, \bibinfo{author}{Francis, B.A.},
  \bibinfo{author}{Tannenbaum, A.R.}, \bibinfo{year}{2013}.
\newblock \bibinfo{title}{Feedback Control Theory}.
\newblock \bibinfo{publisher}{Courier Corporation}.
\bibitem[{Duan et~al.(2024)Duan, Cao, Zheng and Zhao}]{Duan2024LQcontrol}
\bibinfo{author}{Duan, J.}, \bibinfo{author}{Cao, W.}, \bibinfo{author}{Zheng,
  Y.}, \bibinfo{author}{Zhao, L.}, \bibinfo{year}{2024}.
\newblock \bibinfo{title}{On the optimization landscape of dynamic output
  feedback linear quadratic control}.
\newblock \bibinfo{journal}{IEEE Transactions on Automatic Control}
  \bibinfo{volume}{69}, \bibinfo{pages}{920--935}.
\newblock \DOIprefix\doi{10.1109/TAC.2023.3275732}.
\bibitem[{Fam and Meditch(1978)}]{fam1978canonical}
\bibinfo{author}{Fam, A.}, \bibinfo{author}{Meditch, J.}, \bibinfo{year}{1978}.
\newblock \bibinfo{title}{A canonical parameter space for linear systems
  design}.
\newblock \bibinfo{journal}{IEEE Transactions on Automatic Control}
  \bibinfo{volume}{23}, \bibinfo{pages}{454--458}.
\bibitem[{Fazel et~al.(2018a)Fazel, Ge, Kakade and Mesbahi}]{fazel2018global}
\bibinfo{author}{Fazel, M.}, \bibinfo{author}{Ge, R.}, \bibinfo{author}{Kakade,
  S.}, \bibinfo{author}{Mesbahi, M.}, \bibinfo{year}{2018}a.
\newblock \bibinfo{title}{Global convergence of policy gradient methods for the
  linear quadratic regulator}, in: \bibinfo{booktitle}{International conference
  on machine learning}, \bibinfo{organization}{PMLR}. pp.
  \bibinfo{pages}{1467--1476}.
\bibitem[{Fazel et~al.(2018b)Fazel, Ge, Kakade and Mesbahi}]{Fazel2018-pv}
\bibinfo{author}{Fazel, M.}, \bibinfo{author}{Ge, R.}, \bibinfo{author}{Kakade,
  S.M.}, \bibinfo{author}{Mesbahi, M.}, \bibinfo{year}{2018}b.
\newblock \bibinfo{title}{Global convergence of policy gradient methods for the
  linear quadratic regulator}, in: \bibinfo{booktitle}{Proceedings of the 35th
  International Conference on Machine Learning}.
\bibitem[{Feng and Lavaei(2019)}]{feng_exponential_2019}
\bibinfo{author}{Feng, H.}, \bibinfo{author}{Lavaei, J.}, \bibinfo{year}{2019}.
\newblock \bibinfo{title}{On the {Exponential} {Number} of {Connected}
  {Components} for the {Feasible} {Set} of {Optimal} {Decentralized} {Control}
  {Problems}}, in: \bibinfo{booktitle}{2019 {American} {Control} {Conference}
  ({ACC})}, pp. \bibinfo{pages}{1430--1437}.
\newblock \DOIprefix\doi{10.23919/ACC.2019.8814952}.
\bibitem[{Flaxman et~al.(2004)Flaxman, Kalai and McMahan}]{flaxman_online_2004}
\bibinfo{author}{Flaxman, A.D.}, \bibinfo{author}{Kalai, A.T.},
  \bibinfo{author}{McMahan, H.B.}, \bibinfo{year}{2004}.
\newblock \bibinfo{title}{Online convex optimization in the bandit setting:
  gradient descent without a gradient}.
\newblock \URLprefix \url{http://arxiv.org/abs/cs/0408007}.
  \bibinfo{note}{arXiv:cs/0408007}.
\bibitem[{Grondman et~al.()Grondman, Busoniu, Lopes and
  Babusˇka}]{grondman_survey_nodate}
\bibinfo{author}{Grondman, I.}, \bibinfo{author}{Busoniu, L.},
  \bibinfo{author}{Lopes, G.A.D.}, \bibinfo{author}{Babusˇka, R.}, .
\newblock \bibinfo{title}{A {Survey} of {Actor}-{Critic} {Reinforcement}
  {Learning}: {Standard} and {Natural} {Policy} {Gradients}} .
\bibitem[{Guo and Hu(2022)}]{guo2022global}
\bibinfo{author}{Guo, X.}, \bibinfo{author}{Hu, B.}, \bibinfo{year}{2022}.
\newblock \bibinfo{title}{Global convergence of direct policy search for
  state-feedback $\mathcal{H}_\infty$ robust control: A revisit of nonsmooth
  synthesis with {Goldstein} subdifferential}, in: \bibinfo{booktitle}{Advances
  in Neural Information Processing Systems}, \bibinfo{publisher}{Curran
  Associates, Inc.}. pp. \bibinfo{pages}{32801--32815}.
\bibitem[{Guo et~al.(2024)Guo, Keivan, Dullerud, Seiler and
  Hu}]{guo2024complexity}
\bibinfo{author}{Guo, X.}, \bibinfo{author}{Keivan, D.},
  \bibinfo{author}{Dullerud, G.}, \bibinfo{author}{Seiler, P.},
  \bibinfo{author}{Hu, B.}, \bibinfo{year}{2024}.
\newblock \bibinfo{title}{Complexity of derivative-free policy optimization for
  structured $\mathcal{H}_\infty$ control}.
\newblock \bibinfo{journal}{Advances in Neural Information Processing Systems}
  \bibinfo{volume}{36}.
\bibitem[{Hazewinkel(1976)}]{hazewinkel_moduli_1976}
\bibinfo{author}{Hazewinkel, M.}, \bibinfo{year}{1976}.
\newblock \bibinfo{title}{Moduli and canonical forms for linear dynamical
  systems {II}: {The} topological case}.
\newblock \bibinfo{journal}{Mathematical Systems Theory} \bibinfo{volume}{10},
  \bibinfo{pages}{363--385}.
\newblock \URLprefix \url{http://link.springer.com/10.1007/BF01683285},
  \DOIprefix\doi{10.1007/BF01683285}.
\bibitem[{Hazewinkel(1980)}]{byrnes_fine_1980}
\bibinfo{author}{Hazewinkel, M.}, \bibinfo{year}{1980}.
\newblock \bibinfo{title}{({Fine}) {Moduli} ({Spaces}) for {Linear} {Systems}:
  {What} are they and what are they {Good} for?}, in: \bibinfo{editor}{Byrnes,
  C.I.}, \bibinfo{editor}{Martin, C.F.} (Eds.), \bibinfo{booktitle}{Geometrical
  {Methods} for the {Theory} of {Linear} {Systems}}.
  \bibinfo{publisher}{Springer Netherlands}, \bibinfo{address}{Dordrecht}, pp.
  \bibinfo{pages}{125--193}.
\newblock \URLprefix
  \url{http://link.springer.com/10.1007/978-94-009-9082-1_3},
  \DOIprefix\doi{10.1007/978-94-009-9082-1_3}.
\bibitem[{Hazewinkel and Kalman(1976)}]{beckmann_invariants_1976}
\bibinfo{author}{Hazewinkel, M.}, \bibinfo{author}{Kalman, R.E.},
  \bibinfo{year}{1976}.
\newblock \bibinfo{title}{On {Invariants}, {Canonical} {Forms} and {Moduli} for
  {Linear}, {Constant}, {Finite} {Dimensional}, {Dynamical} {Systems}}, in:
  \bibinfo{editor}{Beckmann, M.}, \bibinfo{editor}{Künzi, H.P.},
  \bibinfo{editor}{Marchesini, G.}, \bibinfo{editor}{Mitter, S.K.} (Eds.),
  \bibinfo{booktitle}{Mathematical {Systems} {Theory}}.
  \bibinfo{publisher}{Springer Berlin Heidelberg}, \bibinfo{address}{Berlin,
  Heidelberg}. volume \bibinfo{volume}{131}, pp. \bibinfo{pages}{48--60}.
\newblock \URLprefix
  \url{http://link.springer.com/10.1007/978-3-642-48895-5_4},
  \DOIprefix\doi{10.1007/978-3-642-48895-5_4}.
\bibitem[{Hewer(1971)}]{hewer_iterative_1971}
\bibinfo{author}{Hewer, G.}, \bibinfo{year}{1971}.
\newblock \bibinfo{title}{An iterative technique for the computation of the
  steady state gains for the discrete optimal regulator}.
\newblock \bibinfo{journal}{IEEE Transactions on Automatic Control}
  \bibinfo{volume}{16}, \bibinfo{pages}{382--384}.
\newblock \URLprefix \url{https://ieeexplore.ieee.org/document/1099755},
  \DOIprefix\doi{10.1109/TAC.1971.1099755}. \bibinfo{note}{conference Name:
  IEEE Transactions on Automatic Control}.
\bibitem[{Hu et~al.(2023)Hu, Zhang, Li, Mesbahi, Fazel and
  Ba{\c{s}}ar}]{hu2023toward}
\bibinfo{author}{Hu, B.}, \bibinfo{author}{Zhang, K.}, \bibinfo{author}{Li,
  N.}, \bibinfo{author}{Mesbahi, M.}, \bibinfo{author}{Fazel, M.},
  \bibinfo{author}{Ba{\c{s}}ar, T.}, \bibinfo{year}{2023}.
\newblock \bibinfo{title}{Toward a theoretical foundation of policy
  optimization for learning control policies}.
\newblock \bibinfo{journal}{Annual Review of Control, Robotics, and Autonomous
  Systems} \bibinfo{volume}{6}, \bibinfo{pages}{123--158}.
\bibitem[{Hu and Zheng(2022)}]{hu2022connectivity}
\bibinfo{author}{Hu, B.}, \bibinfo{author}{Zheng, Y.}, \bibinfo{year}{2022}.
\newblock \bibinfo{title}{Connectivity of the feasible and sublevel sets of
  dynamic output feedback control with robustness constraints}.
\newblock \bibinfo{journal}{IEEE Control Systems Letters} \bibinfo{volume}{7},
  \bibinfo{pages}{442--447}.
\bibitem[{Hyland and Bernstein(1984)}]{hyland1984optimal}
\bibinfo{author}{Hyland, D.}, \bibinfo{author}{Bernstein, D.},
  \bibinfo{year}{1984}.
\newblock \bibinfo{title}{The optimal projection equations for fixed-order
  dynamic compensation}.
\newblock \bibinfo{journal}{IEEE Transactions on Automatic Control}
  \bibinfo{volume}{29}, \bibinfo{pages}{1034--1037}.
\bibitem[{Isidori(2013)}]{Isidori2013-dj}
\bibinfo{author}{Isidori, A.}, \bibinfo{year}{2013}.
\newblock \bibinfo{title}{Nonlinear Control Systems}.
\newblock \bibinfo{publisher}{Springer Science \& Business Media}.
\bibitem[{Jin et~al.(2017)Jin, Ge, Netrapalli, Kakade and
  Jordan}]{jin2017escape}
\bibinfo{author}{Jin, C.}, \bibinfo{author}{Ge, R.},
  \bibinfo{author}{Netrapalli, P.}, \bibinfo{author}{Kakade, S.M.},
  \bibinfo{author}{Jordan, M.I.}, \bibinfo{year}{2017}.
\newblock \bibinfo{title}{How to escape saddle points efficiently}, in:
  \bibinfo{booktitle}{International conference on machine learning},
  \bibinfo{organization}{PMLR}. pp. \bibinfo{pages}{1724--1732}.
\bibitem[{Jury(1964)}]{jury1964theory}
\bibinfo{author}{Jury, E.I.}, \bibinfo{year}{1964}.
\newblock \bibinfo{title}{Theory and application of the z-transform method} .
\bibitem[{Kakade(2001)}]{kakade_natural_2001}
\bibinfo{author}{Kakade, S.M.}, \bibinfo{year}{2001}.
\newblock \bibinfo{title}{A {Natural} {Policy} {Gradient}}, in:
  \bibinfo{booktitle}{Advances in {Neural} {Information} {Processing}
  {Systems}}, \bibinfo{publisher}{MIT Press}.
\newblock \URLprefix
  \url{https://proceedings.neurips.cc/paper_files/paper/2001/hash/4b86abe48d358ecf194c56c69108433e-Abstract.html}.
\bibitem[{Kakade(2002)}]{Kakade2002-fc}
\bibinfo{author}{Kakade, S.M.}, \bibinfo{year}{2002}.
\newblock \bibinfo{title}{A natural policy gradient}, in:
  \bibinfo{editor}{Dietterich, T.G.}, \bibinfo{editor}{Becker, S.},
  \bibinfo{editor}{Ghahramani, Z.} (Eds.), \bibinfo{booktitle}{Advances in
  Neural Information Processing Systems 14}. \bibinfo{publisher}{MIT Press},
  pp. \bibinfo{pages}{1531--1538}.
\bibitem[{Keivan et~al.(2022)Keivan, Havens, Seiler, Dullerud and
  Hu}]{keivan2022model}
\bibinfo{author}{Keivan, D.}, \bibinfo{author}{Havens, A.},
  \bibinfo{author}{Seiler, P.}, \bibinfo{author}{Dullerud, G.},
  \bibinfo{author}{Hu, B.}, \bibinfo{year}{2022}.
\newblock \bibinfo{title}{Model-free $\mu$ synthesis via adversarial
  reinforcement learning}, in: \bibinfo{booktitle}{2022 American Control
  Conference (ACC)}, \bibinfo{organization}{IEEE}. pp.
  \bibinfo{pages}{3335--3341}.
\bibitem[{Kraisler and Mesbahi(2024)}]{kraisler2024output}
\bibinfo{author}{Kraisler, S.}, \bibinfo{author}{Mesbahi, M.},
  \bibinfo{year}{2024}.
\newblock \bibinfo{title}{Output-feedback synthesis orbit geometry: Quotient
  manifolds and lqg direct policy optimization}.
\newblock \bibinfo{journal}{arXiv preprint arXiv:2403.17157} .
\bibitem[{Krishnaprasad and Martin(1983)}]{krishnaprasad_families_1983a}
\bibinfo{author}{Krishnaprasad, P.S.}, \bibinfo{author}{Martin, C.F.},
  \bibinfo{year}{1983}.
\newblock \bibinfo{title}{On families of systems and deformations}.
\newblock \bibinfo{journal}{International Journal of Control}
  \bibinfo{volume}{38}, \bibinfo{pages}{1055--1079}.
\bibitem[{Lee(2010)}]{lee_introduction_2010}
\bibinfo{author}{Lee, J.}, \bibinfo{year}{2010}.
\newblock \bibinfo{title}{Introduction to {Topological} {Manifolds}}.
\newblock \bibinfo{edition}{2nd} ed., \bibinfo{publisher}{Springer Science \&
  Business Media}.
\bibitem[{Lee(2018)}]{lee_introduction_2018}
\bibinfo{author}{Lee, J.M.}, \bibinfo{year}{2018}.
\newblock \bibinfo{title}{Introduction to {Riemannian} {Manifolds}}.
\newblock \bibinfo{edition}{2nd} ed., \bibinfo{publisher}{Cham, Switzerland:
  Springer Nature}.
\bibitem[{Lewis(2007)}]{lewis2007nonsmooth}
\bibinfo{author}{Lewis, A.S.}, \bibinfo{year}{2007}.
\newblock \bibinfo{title}{Nonsmooth optimization and robust control}.
\newblock \bibinfo{journal}{Annual Reviews in Control} \bibinfo{volume}{31},
  \bibinfo{pages}{167--177}.
\bibitem[{Liberzon(2011)}]{Liberzon2011-en}
\bibinfo{author}{Liberzon, D.}, \bibinfo{year}{2011}.
\newblock \bibinfo{title}{Calculus of Variations and Optimal Control Theory: A
  Concise Introduction}.
\newblock \bibinfo{publisher}{Princeton University Press}.
\bibitem[{Malik et~al.(2019)Malik, Pananjady, Bhatia, Khamaru, Bartlett and
  Wainwright}]{Malik2019-bs}
\bibinfo{author}{Malik, D.}, \bibinfo{author}{Pananjady, A.},
  \bibinfo{author}{Bhatia, K.}, \bibinfo{author}{Khamaru, K.},
  \bibinfo{author}{Bartlett, P.}, \bibinfo{author}{Wainwright, M.},
  \bibinfo{year}{2019}.
\newblock \bibinfo{title}{{Derivative-Free} methods for policy optimization:
  Guarantees for linear quadratic systems}, in: \bibinfo{editor}{Chaudhuri,
  K.}, \bibinfo{editor}{Sugiyama, M.} (Eds.), \bibinfo{booktitle}{Proceedings
  of the {Twenty-Second} International Conference on Artificial Intelligence
  and Statistics}, \bibinfo{publisher}{PMLR}. pp. \bibinfo{pages}{2916--2925}.
\bibitem[{Mohammadi et~al.(2020)Mohammadi, Soltanolkotabi and
  Jovanovic}]{Mohammadi2020-iv}
\bibinfo{author}{Mohammadi, H.}, \bibinfo{author}{Soltanolkotabi, M.},
  \bibinfo{author}{Jovanovic, M.R.}, \bibinfo{year}{2020}.
\newblock \bibinfo{title}{Random search for learning the linear quadratic
  regulator}, in: \bibinfo{booktitle}{2020 American Control Conference
  ({ACC})}, \bibinfo{publisher}{IEEE}. pp. \bibinfo{pages}{4798--4803}.
\bibitem[{Mohammadi et~al.(2021a)Mohammadi, Soltanolkotabi and
  Jovanovi{\'c}}]{Mohammadi2021-vs}
\bibinfo{author}{Mohammadi, H.}, \bibinfo{author}{Soltanolkotabi, M.},
  \bibinfo{author}{Jovanovi{\'c}, M.R.}, \bibinfo{year}{2021}a.
\newblock \bibinfo{title}{{Model-Free} linear quadratic regulator}, in:
  \bibinfo{editor}{Vamvoudakis, K.G.}, \bibinfo{editor}{Wan, Y.},
  \bibinfo{editor}{Lewis, F.L.}, \bibinfo{editor}{Cansever, D.} (Eds.),
  \bibinfo{booktitle}{Handbook of Reinforcement Learning and Control}.
  \bibinfo{publisher}{Springer International Publishing},
  \bibinfo{address}{Cham}, pp. \bibinfo{pages}{173--185}.
\bibitem[{Mohammadi et~al.(2021b)Mohammadi, Soltanolkotabi and
  Jovanovi{\'c}}]{mohammadi2021lack}
\bibinfo{author}{Mohammadi, H.}, \bibinfo{author}{Soltanolkotabi, M.},
  \bibinfo{author}{Jovanovi{\'c}, M.R.}, \bibinfo{year}{2021}b.
\newblock \bibinfo{title}{On the lack of gradient domination for linear
  quadratic gaussian problems with incomplete state information}, in:
  \bibinfo{booktitle}{2021 60th IEEE Conference on Decision and Control (CDC)},
  \bibinfo{organization}{IEEE}. pp. \bibinfo{pages}{1120--1124}.
\bibitem[{Mohammadi et~al.(2021c)Mohammadi, Zare, Soltanolkotabi and
  Jovanovi{\'c}}]{mohammadi2021convergence}
\bibinfo{author}{Mohammadi, H.}, \bibinfo{author}{Zare, A.},
  \bibinfo{author}{Soltanolkotabi, M.}, \bibinfo{author}{Jovanovi{\'c}, M.R.},
  \bibinfo{year}{2021}c.
\newblock \bibinfo{title}{Convergence and sample complexity of gradient methods
  for the model-free linear--quadratic regulator problem}.
\newblock \bibinfo{journal}{IEEE Transactions on Automatic Control}
  \bibinfo{volume}{67}, \bibinfo{pages}{2435--2450}.
\bibitem[{Nijmeijer and Schaft(1990)}]{Nijmeijer1990-xr}
\bibinfo{author}{Nijmeijer, H.}, \bibinfo{author}{Schaft, A.},
  \bibinfo{year}{1990}.
\newblock \bibinfo{title}{Nonlinear Dynamical Control Systems}.
\newblock \bibinfo{publisher}{Springer New York}.
\bibitem[{Ober(1987)}]{ober1987topology}
\bibinfo{author}{Ober, R.J.}, \bibinfo{year}{1987}.
\newblock \bibinfo{title}{Topology of the set of asymptotically stable minimal
  systems}.
\newblock \bibinfo{journal}{International Journal of Control}
  \bibinfo{volume}{46}, \bibinfo{pages}{263--280}.
\bibitem[{Ohara and Amari(1992a)}]{ohara1992differential}
\bibinfo{author}{Ohara, A.}, \bibinfo{author}{Amari, S.i.},
  \bibinfo{year}{1992}a.
\newblock \bibinfo{title}{Differential geometric structures of stable state
  feedback systems with dual connections}.
\newblock \bibinfo{journal}{IFAC Proceedings Volumes} \bibinfo{volume}{25},
  \bibinfo{pages}{176--179}.
\bibitem[{Ohara and Amari(1992b)}]{ohara_differential_1992}
\bibinfo{author}{Ohara, A.}, \bibinfo{author}{Amari, S.I.},
  \bibinfo{year}{1992}b.
\newblock \bibinfo{title}{Differential {Geometric} {Structures} of {Stable}
  {State} {Feedback} {Systems} with {Dual} {Connections}}.
\newblock \bibinfo{journal}{IFAC Proceedings} \bibinfo{volume}{25},
  \bibinfo{pages}{176--179}.
\newblock \URLprefix
  \url{https://www.sciencedirect.com/science/article/pii/S147466701749745X},
  \DOIprefix\doi{https://doi.org/10.1016/S1474-6670(17)49745-X}.
\bibitem[{Parks(1962)}]{parks1962new}
\bibinfo{author}{Parks, P.C.}, \bibinfo{year}{1962}.
\newblock \bibinfo{title}{A new proof of the {Routh-Hurwitz} stability
  criterion using the second method of liapunov}, in:
  \bibinfo{booktitle}{Mathematical Proceedings of the Cambridge Philosophical
  Society}, \bibinfo{organization}{Cambridge University Press}. pp.
  \bibinfo{pages}{694--702}.
\bibitem[{Rantzer(1996)}]{rantzer1996kalman}
\bibinfo{author}{Rantzer, A.}, \bibinfo{year}{1996}.
\newblock \bibinfo{title}{On the {Kalman--Yakubovich--Popov} lemma}.
\newblock \bibinfo{journal}{Systems \& Control Letters} \bibinfo{volume}{28},
  \bibinfo{pages}{7--10}.
\bibitem[{Recht(2019)}]{recht2019tour}
\bibinfo{author}{Recht, B.}, \bibinfo{year}{2019}.
\newblock \bibinfo{title}{A tour of reinforcement learning: The view from
  continuous control}.
\newblock \bibinfo{journal}{Annual Review of Control, Robotics, and Autonomous
  Systems} \bibinfo{volume}{2}, \bibinfo{pages}{253--279}.
\bibitem[{Schulman et~al.(2017)Schulman, Wolski, Dhariwal, Radford and
  Klimov}]{schulman_proximal_2017}
\bibinfo{author}{Schulman, J.}, \bibinfo{author}{Wolski, F.},
  \bibinfo{author}{Dhariwal, P.}, \bibinfo{author}{Radford, A.},
  \bibinfo{author}{Klimov, O.}, \bibinfo{year}{2017}.
\newblock \bibinfo{title}{Proximal {Policy} {Optimization} {Algorithms}}.
\newblock \URLprefix \url{http://arxiv.org/abs/1707.06347}.
  \bibinfo{note}{arXiv:1707.06347 [cs]}.
\bibitem[{Skogestad and Postlethwaite(2005)}]{Skogestad2005-wi}
\bibinfo{author}{Skogestad, S.}, \bibinfo{author}{Postlethwaite, I.},
  \bibinfo{year}{2005}.
\newblock \bibinfo{title}{Multivariable Feedback Control: Analysis and Design}.
\newblock \bibinfo{publisher}{John Wiley \& Sons}.
\bibitem[{Sontag(2013)}]{Sontag2013-hk}
\bibinfo{author}{Sontag, E.D.}, \bibinfo{year}{2013}.
\newblock \bibinfo{title}{Mathematical Control Theory: Deterministic Finite
  Dimensional Systems}.
\newblock \bibinfo{publisher}{Springer Science \& Business Media}.
\bibitem[{Spall(1998)}]{spall_overview_1998}
\bibinfo{author}{Spall, J.C.}, \bibinfo{year}{1998}.
\newblock \bibinfo{title}{An {Overview} of the {Simultaneous} {Perturbation}
  {Method}}.
\newblock \bibinfo{journal}{JOHNS HOPKINS APL TECHNICAL DIGEST}
  \bibinfo{volume}{19}.
\bibitem[{Sutton et~al.(1999)Sutton, McAllester, Singh and
  Mansour}]{sutton_policy_1999}
\bibinfo{author}{Sutton, R.S.}, \bibinfo{author}{McAllester, D.},
  \bibinfo{author}{Singh, S.}, \bibinfo{author}{Mansour, Y.},
  \bibinfo{year}{1999}.
\newblock \bibinfo{title}{Policy {Gradient} {Methods} for {Reinforcement}
  {Learning} with {Function} {Approximation}}, in: \bibinfo{booktitle}{Advances
  in {Neural} {Information} {Processing} {Systems}}, \bibinfo{publisher}{MIT
  Press}.
\newblock \URLprefix
  \url{https://proceedings.neurips.cc/paper/1999/hash/464d828b85b0bed98e80ade0a5c43b0f-Abstract.html}.
\bibitem[{Takakura and Sato(2024)}]{takakura_structured_2024}
\bibinfo{author}{Takakura, S.}, \bibinfo{author}{Sato, K.},
  \bibinfo{year}{2024}.
\newblock \bibinfo{title}{Structured {Output} {Feedback} {Control} for {Linear}
  {Quadratic} {Regulator} {Using} {Policy} {Gradient} {Method}}.
\newblock \bibinfo{journal}{IEEE Transactions on Automatic Control}
  \bibinfo{volume}{69}, \bibinfo{pages}{363--370}.
\newblock \URLprefix \url{https://ieeexplore.ieee.org/document/10091214/},
  \DOIprefix\doi{10.1109/TAC.2023.3264176}.
\bibitem[{Talebi et~al.(2022)Talebi, Alemzadeh, Rahimi and
  Mesbahi}]{talebi_regularizability_2022}
\bibinfo{author}{Talebi, S.}, \bibinfo{author}{Alemzadeh, S.},
  \bibinfo{author}{Rahimi, N.}, \bibinfo{author}{Mesbahi, M.},
  \bibinfo{year}{2022}.
\newblock \bibinfo{title}{On {Regularizability} and {Its} {Application} to
  {Online} {Control} of {Unstable} {LTI} {Systems}}.
\newblock \bibinfo{journal}{IEEE Transactions on Automatic Control}
  \bibinfo{volume}{67}, \bibinfo{pages}{6413--6428}.
\newblock \URLprefix
  \url{https://ieeexplore.ieee.org/abstract/document/9627635},
  \DOIprefix\doi{10.1109/TAC.2021.3131148}. \bibinfo{note}{conference Name:
  IEEE Transactions on Automatic Control}.
\bibitem[{Talebi and Mesbahi(2022)}]{talebi_riemannian_2022}
\bibinfo{author}{Talebi, S.}, \bibinfo{author}{Mesbahi, M.},
  \bibinfo{year}{2022}.
\newblock \bibinfo{title}{Riemannian {Constrained} {Policy} {Optimization} via
  {Geometric} {Stability} {Certificates}}, in: \bibinfo{booktitle}{2022 {IEEE}
  61st {Conference} on {Decision} and {Control} ({CDC})}, pp.
  \bibinfo{pages}{1472--1478}.
\newblock \DOIprefix\doi{10.1109/CDC51059.2022.9992877}.
\bibitem[{Talebi and Mesbahi(2023)}]{talebi_policy_2023}
\bibinfo{author}{Talebi, S.}, \bibinfo{author}{Mesbahi, M.},
  \bibinfo{year}{2023}.
\newblock \bibinfo{title}{Policy {Optimization} over {Submanifolds} for
  {Linearly} {Constrained} {Feedback} {Synthesis}}.
\newblock \bibinfo{journal}{IEEE Transactions on Automatic Control} ,
  \bibinfo{pages}{1--16}\URLprefix
  \url{https://ieeexplore.ieee.org/document/10224332},
  \DOIprefix\doi{10.1109/TAC.2023.3306384}. \bibinfo{note}{conference Name:
  IEEE Transactions on Automatic Control}.
\bibitem[{Talebi et~al.(2023)Talebi, Taghvaei and
  Mesbahi}]{talebi_data-driven_2023}
\bibinfo{author}{Talebi, S.}, \bibinfo{author}{Taghvaei, A.},
  \bibinfo{author}{Mesbahi, M.}, \bibinfo{year}{2023}.
\newblock \bibinfo{title}{Data-driven {Optimal} {Filtering} for {Linear}
  {Systems} with {Unknown} {Noise} {Covariances}}.
\newblock \bibinfo{journal}{Advances in Neural Information Processing Systems}
  \bibinfo{volume}{36}, \bibinfo{pages}{69546--69585}.
\newblock \URLprefix
  \url{https://proceedings.neurips.cc/paper_files/paper/2023/hash/dbe8185809cb7032ec7ec6e365e3ed3b-Abstract-Conference.html}.
\bibitem[{Tang and Zheng(2023)}]{tang2023global}
\bibinfo{author}{Tang, Y.}, \bibinfo{author}{Zheng, Y.}, \bibinfo{year}{2023}.
\newblock \bibinfo{title}{On the global optimality of direct policy search for
  nonsmooth $\mathcal{H}_\infty$ output-feedback control}, in:
  \bibinfo{booktitle}{2023 62nd IEEE Conference on Decision and Control (CDC)},
  \bibinfo{organization}{IEEE}. pp. \bibinfo{pages}{6148--6153}.
\bibitem[{Tang et~al.(2023)Tang, Zheng and Li}]{tang2023analysis}
\bibinfo{author}{Tang, Y.}, \bibinfo{author}{Zheng, Y.}, \bibinfo{author}{Li,
  N.}, \bibinfo{year}{2023}.
\newblock \bibinfo{title}{Analysis of the optimization landscape of linear
  quadratic gaussian (lqg) control}.
\newblock \bibinfo{journal}{Mathematical Programming} \bibinfo{volume}{202},
  \bibinfo{pages}{399--444}.
\bibitem[{Tannenbaum(1981)}]{tannenbaum_invariance_1981}
\bibinfo{author}{Tannenbaum, A.R.}, \bibinfo{year}{1981}.
\newblock \bibinfo{title}{Invariance and system theory: algebraic and geometric
  aspects}.
\newblock Number \bibinfo{number}{845} in \bibinfo{series}{Lecture notes in
  mathematics}, \bibinfo{publisher}{Springer}, \bibinfo{address}{Berlin}.
\bibitem[{Trentelman et~al.(2001)Trentelman, Stoorvogel and
  Hautus}]{Trentelman2001-ck}
\bibinfo{author}{Trentelman, H.L.}, \bibinfo{author}{Stoorvogel, A.A.},
  \bibinfo{author}{Hautus, M.}, \bibinfo{year}{2001}.
\newblock \bibinfo{title}{Control theory for linear systems}.
\bibitem[{Umenberger et~al.(2022)Umenberger, Simchowitz, Perdomo, Zhang and
  Tedrake}]{umenberger2022globally}
\bibinfo{author}{Umenberger, J.}, \bibinfo{author}{Simchowitz, M.},
  \bibinfo{author}{Perdomo, J.}, \bibinfo{author}{Zhang, K.},
  \bibinfo{author}{Tedrake, R.}, \bibinfo{year}{2022}.
\newblock \bibinfo{title}{Globally convergent policy search for output
  estimation}.
\newblock \bibinfo{journal}{Advances in Neural Information Processing Systems}
  \bibinfo{volume}{35}, \bibinfo{pages}{22778--22790}.
\bibitem[{Williams(1992)}]{williams_simple_1992}
\bibinfo{author}{Williams, R.J.}, \bibinfo{year}{1992}.
\newblock \bibinfo{title}{Simple statistical gradient-following algorithms for
  connectionist reinforcement learning}.
\newblock \bibinfo{journal}{Machine Learning} \bibinfo{volume}{8},
  \bibinfo{pages}{229--256}.
\newblock \URLprefix \url{https://doi.org/10.1007/BF00992696},
  \DOIprefix\doi{10.1007/BF00992696}.
\bibitem[{Wonham(1985)}]{wonham_w_murray_linear_1985}
\bibinfo{author}{Wonham, W.M.}, \bibinfo{year}{1985}.
\newblock \bibinfo{title}{Linear {Multivariable} {Control}}.
\newblock \bibinfo{publisher}{Springer Science+Business Media, LLC},
  \bibinfo{address}{New York}.
\newblock \URLprefix
  \url{https://link-springer-com.ezp-prod1.hul.harvard.edu/book/10.1007/978-1-4612-1082-5}.
\bibitem[{Yu et~al.(2023)Yu, Gupta and Wierman}]{yu_online_2023}
\bibinfo{author}{Yu, J.}, \bibinfo{author}{Gupta, V.},
  \bibinfo{author}{Wierman, A.}, \bibinfo{year}{2023}.
\newblock \bibinfo{title}{Online {Stabilization} of {Unknown} {Linear}
  {Time}-{Varying} {Systems}}.
\newblock \URLprefix \url{https://arxiv.org/abs/2304.02878v2}.
\bibitem[{Zabczyk(2020)}]{Zabczyk_undated-jw}
\bibinfo{author}{Zabczyk, J.}, \bibinfo{year}{2020}.
\newblock \bibinfo{title}{Mathematical Control Theory}.
\newblock \bibinfo{publisher}{Birkh{\"a}user}.
\bibitem[{Zhang et~al.(2020)Zhang, Hu and Basar}]{zhang2020policy}
\bibinfo{author}{Zhang, K.}, \bibinfo{author}{Hu, B.}, \bibinfo{author}{Basar,
  T.}, \bibinfo{year}{2020}.
\newblock \bibinfo{title}{Policy optimization for $\mathcal{H}_2$ linear
  control with $\mathcal{H}_\infty$ robustness guarantee: Implicit
  regularization and global convergence}, in: \bibinfo{booktitle}{Learning for
  Dynamics and Control}, \bibinfo{organization}{PMLR}. pp.
  \bibinfo{pages}{179--190}.
\bibitem[{Zhang et~al.(2023)Zhang, Hu and Ba{\c{s}}ar}]{zhang2023learning}
\bibinfo{author}{Zhang, X.}, \bibinfo{author}{Hu, B.},
  \bibinfo{author}{Ba{\c{s}}ar, T.}, \bibinfo{year}{2023}.
\newblock \bibinfo{title}{Learning the kalman filter with fine-grained sample
  complexity}, in: \bibinfo{booktitle}{2023 American Control Conference (ACC)},
  \bibinfo{organization}{IEEE}. pp. \bibinfo{pages}{4549--4554}.
\bibitem[{Zhao et~al.(2023)Zhao, You and Ba{\c{s}}ar}]{zhao2023global}
\bibinfo{author}{Zhao, F.}, \bibinfo{author}{You, K.},
  \bibinfo{author}{Ba{\c{s}}ar, T.}, \bibinfo{year}{2023}.
\newblock \bibinfo{title}{Global convergence of policy gradient primal-dual
  methods for risk-constrained lqrs}.
\newblock \bibinfo{journal}{IEEE Transactions on Automatic Control} .
\bibitem[{Zheng et~al.(2023)Zheng, Pai and Tang}]{zheng2023benign}
\bibinfo{author}{Zheng, Y.}, \bibinfo{author}{Pai, C.f.},
  \bibinfo{author}{Tang, Y.}, \bibinfo{year}{2023}.
\newblock \bibinfo{title}{Benign nonconvex landscapes in optimal and robust
  control, {Part I}: Global optimality}.
\newblock \bibinfo{journal}{arXiv preprint arXiv:2312.15332} .
\bibitem[{Zheng et~al.(2024)Zheng, Pai and Tang}]{zheng2024benign}
\bibinfo{author}{Zheng, Y.}, \bibinfo{author}{Pai, C.f.},
  \bibinfo{author}{Tang, Y.}, \bibinfo{year}{2024}.
\newblock \bibinfo{title}{Benign nonconvex landscapes in optimal and robust
  control, {Part II}: Extended convex lifting}.
\newblock \bibinfo{journal}{preprint} .
\bibitem[{Zheng et~al.(2022)Zheng, Sun, Fazel and Li}]{zheng2022escaping}
\bibinfo{author}{Zheng, Y.}, \bibinfo{author}{Sun, Y.}, \bibinfo{author}{Fazel,
  M.}, \bibinfo{author}{Li, N.}, \bibinfo{year}{2022}.
\newblock \bibinfo{title}{Escaping high-order saddles in policy optimization
  for linear quadratic gaussian (lqg) control}, in: \bibinfo{booktitle}{2022
  IEEE 61st Conference on Decision and Control (CDC)},
  \bibinfo{organization}{IEEE}. pp. \bibinfo{pages}{5329--5334}.
\bibitem[{Zhou et~al.(1996)Zhou, Doyle and Glover}]{zhou1996robust}
\bibinfo{author}{Zhou, K.}, \bibinfo{author}{Doyle, J.C.},
  \bibinfo{author}{Glover, K.}, \bibinfo{year}{1996}.
\newblock \bibinfo{title}{Robust and Optimal Control}.
\newblock \bibinfo{publisher}{Prentice Hall}.
\bibitem[{Zhou and Khargonekar(1988)}]{zhou1988algebraic}
\bibinfo{author}{Zhou, K.}, \bibinfo{author}{Khargonekar, P.P.},
  \bibinfo{year}{1988}.
\newblock \bibinfo{title}{An algebraic riccati equation approach to
  $\mathcal{H}_\infty$ optimization}.
\newblock \bibinfo{journal}{Systems \& Control Letters} \bibinfo{volume}{11},
  \bibinfo{pages}{85--91}.

\end{thebibliography}
\end{document}